\documentclass[a4paper]{amsart}
\usepackage{times}
\usepackage{amsfonts,amscd,amsmath,amssymb,amsthm}
\usepackage[bbgreekl]{mathbbol}

\numberwithin{equation}{section}

\newcommand{\sk}{\mathsf{k}}

\renewcommand{\Im}{{\ensuremath{\mathrm{Im\,}}}}
\renewcommand{\Re}{{\ensuremath{\mathrm{Re\,}}}}

\DeclareSymbolFont{SY}{U}{psy}{m}{n}
\DeclareMathSymbol{\emptyset}{\mathord}{SY}{'306}

\DeclareMathOperator{\Ran}{Ran} \DeclareMathOperator{\Ker}{Ker}
 \DeclareMathOperator{\spec}{spec}
 \DeclareMathOperator{\Dom}{Dom}
\DeclareMathOperator{\sign}{sign}

\DeclareMathOperator{\ind}{ind}

\DeclareMathSymbol{\newtimes}{\mathbin}{SY}{'264}
\DeclareMathOperator*{\Bigtimes}{\newtimes}

\newcommand{\diag}{\mathrm{diag}}

\newcommand{\ul}[1]{\ensuremath\underline{#1}}

\newcommand{\R}{\mathbb{R}}

\newcommand{\C}{\mathbb{C}}

\newcommand{\1}{\mathbb{I}}

\newcommand{\fH}{\mathfrak{H}}

\newcommand{\fS}{\mathfrak{S}}

\newcommand{\cC}{{\mathcal C}}
\newcommand{\cD}{{\mathcal D}}
\newcommand{\cE}{{\mathcal E}}

\newcommand{\cG}{{\mathcal G}}
\newcommand{\cH}{{\mathcal H}}
\newcommand{\cI}{{\mathcal I}}
\newcommand{\cJ}{{\mathcal J}}
\newcommand{\cK}{{\mathcal K}}
\newcommand{\cL}{{\mathcal L}}
\newcommand{\cM}{{\mathcal M}}
\newcommand{\cN}{{\mathcal N}}

\newcommand{\cW}{{\mathcal W}}

\newcommand{\bw}{\mathbf{w}}

\newcommand{\ii}{\mathrm{i}}
\newcommand{\e}{\mathrm{e}}

\newcommand{\CDelta}{\mathbb{\Delta}}

{\bf}{\it}
{\bf}{\it}
{\bf}{\it}
{\bf}{\it}
{\bf}{\it}

\newtheorem{theorem}{Theorem}[section]{\bf}{\it}
\newtheorem{proposition}[theorem]{Proposition}{\bf}{\it}
\newtheorem{corollary}[theorem]{Corollary}{\bf}{\it}
\newtheorem{lemma}[theorem]{Lemma}{\bf}{\it}
\newtheorem{remark}[theorem]{Remark}{\it}{\rm}
\newtheorem{definition}[theorem]{Definition}{\bf}{\it}
\newtheorem{assumption}[theorem]{Assumption}{\bf}{\it}
\theoremstyle{remark}
\newtheorem{example}[theorem]{Example}{\bf}{\rm}

\title[Contraction Semigroups on Metric Graphs]{Contraction Semigroups on Metric Graphs}

\author[V. Kostrykin]{Vadim Kostrykin}
\address{Vadim Kostrykin\\ Institut f\"{u}r Mathematik, Technische Universit\"{a}t Clausthal,
Erz\-stra{\ss}e 1, D-38678 Clausthal-Zellerfeld, Germany;\newline
Fraunhofer-Institut f\"{u}r Laser\-technik, Steinbachstra{\ss}e 15, 52074 Aachen, Germany}
\email{kostrykin@math.tu-clausthal.de, kostrykin@t-online.de}

\author[J. Potthoff]{J\"{u}rgen Potthoff}
\address{J\"{u}rgen Potthoff\\ Institut f\"{u}r Mathematik, Universit\"{a}t Mannheim, D-68131 Mannheim, Germany}
\email{potthoff@math.uni-mannheim.de}

\author[R. Schrader]{Robert Schrader}
\address{Robert Schrader\\ Institut f\"{u}r Theoretische Physik\\
Freie Universit\"{a}t Berlin, Arnim\-allee 14\\ D-14195 Berlin, Germany}
\email{schrader@physik.fu-berlin.de}

\dedicatory{Dedicated to Volker Enss on the occasion of his 65-th birthday}

\copyrightinfo{2008}{V.~Kostrykin, J.~Potthoff, R.~Schrader}

\subjclass[2000]{34B45, 47D06, 47B44}

\keywords{Differential operators on metric graphs,
accretive and dissipative extensions of symmetric operators, contraction semigroups, Feller semigroups.}

\begin{document}

\begin{abstract}
The main objective of the present work is to study contraction semigroups generated by Laplace operators on metric graphs, which are not
necessarily self-adjoint. We prove criteria for such semigroups to be continuity and positivity preserving. Also we provide a characterization of generators of Feller semigroups on metric graphs.
\end{abstract}

\maketitle

\section{Introduction}\label{sec:intro}

Metric graphs or networks are one-dimensional piecewise linear spaces with
singularities at the vertices. Alternatively, a metric graph is a metric
space which can be written as a union of finitely many intervals, which are
either compact or $[0,+\infty)$; any two of these intervals are either
disjoint or intersect only in one or both of their endpoints. It is natural
to call the metric graph compact if all its edges have finite length.

The increasing interest in the theory of differential operators on
metric graphs is motivated mainly by two reasons. The first reason
is that such operators arise in a variety of applications. We
refer the reader to the review \cite{Kuchment:00}, where a number
of models arising in physics, chemistry, and engineering are
discussed, as well as to original works \cite{Badalin:1}, \cite{Badalin:2}, \cite{Kuchment:Post}, where modeling of carbon nano-structures is discussed. References \cite{Cardanobile}, \cite{Carlson}, \cite{Nicaise:0} address signal transmission in biological neural networks and blood flow in the human arterial systems. The second reason is
purely mathematical: It is intriguing to study the interrelation
between the spectra of these operators and topological or
combinatorial properties of the underlying graph. Similar
interrelations are studied in spectral geometry for differential
operators on Riemannian manifolds (see, e.g.~\cite{Chavel},
\cite{Gilkey}, \cite{Gilkey:2}) and in spectral graph theory for difference
operators on combinatorial graphs (see, e.g.~\cite{Chung}, \cite{Verdiere}).

Metric graphs take an intermediate position between manifolds and
combinatorial graphs. References \cite{KPS1}, \cite{Roth:1}, \cite{Roth:2} provide a Selberg-type trace formula for semigroups generated by a class of self-adjoint Laplace operators on metric graphs which establishes a connection between the trace of the semigroup and
cycles on the graph as well as its Euler characteristics (see also \cite{Kurasov}). Index theorems for such semigroups
have been proved in \cite{Fulling}. These results have the well-known analogues in spectral geometry. On the other hand, for compact graphs with equal edge lengths and special boundary conditions at the vertices the spectrum of the differential Laplace operator is the preimage of the combinatorial spectrum under a certain
entire function (see, e.g., \cite{Pankrashkin}). Using this correspondence some results of the spectral graph theory for combinatorial Laplacians can be translated in this case to differential Laplace operators.

In the setting of the Hilbert space theory, semigroups generated by self-adjoint differential operators on metric graphs in special cases have been studied in \cite{Angad-Guar}, \cite{Gaveau:1}, \cite{Gaveau:2}, \cite{Okada}. Reference \cite{KS9} seems to be the first work, where a systematic study of semigroups on metric graphs has been undertaken. In particular, this reference provides criteria for a self-adjoint Laplace operator to generate a contraction and positivity preserving semigroup.

The main objective of the present paper is to study semigroups
generated by general, not necessarily self-adjoint Laplace
operators on metric graphs. There are several motivations to study such semigroups.

First, non-self-adjoint differential operators naturally appear in different models. In particular, initially motivated by neurobiological problems, parabolic equations on (finite or infinite) metric graphs attract research attention since more than 20 years
(see, e.g., \cite{Ali:Mehmeti}, \cite{Below}, \cite{Mugnolo:4}, \cite{Lumer:2}, \cite{Mugnolo:2}, \cite{Mugnolo:3}, \cite{Nicaise:0} and references quoted therein). Here the theory of semigroups on metric graphs plays a central role.

Second, positivity preserving contraction semigroups on the Banach space of continuous functions (that is, Feller semigroups) are related to strongly Markovian stochastic processes on metric graphs. In particular, the transition density of a stochastic processes is given by the integral kernel of the corresponding semigroups. Without attempting to give a complete review of the work on stochastic processes on metric graphs we mention the papers \cite{Baxter}, \cite{Freidlin:1}, \cite{Freidlin:2}, \cite{Walsh}. In our forthcoming article \cite{KPS2} we give a complete characterization and construction of all Brownian motions on metric graphs, that is,
of all path continuous strongly Markovian stochastic
processes which away from the vertices are equivalent
to a standard Brownian motion.

The work is organized as follows. Notation and main results are presented in Section \ref{sec:basic}.  In Section \ref{sec:diss:accr} accretive and dissipative Laplace operators are studied. In particular, we provide a characterization of maximal accretive operators, which by the Lumer-Phillips theorem are infinitesimal generators of contraction semigroups. In Sections \ref{sec:pres:cont} and \ref{sec:pos:preserv} we prove criteria for the boundary conditions at the vertices of the graph ensuring that the contraction semigroup generated by the corresponding Laplace operator is
positivity and continuity preserving. The semigroup theory in the Hilbert space developed here can be used to
study the semigroups on other function spaces on metric graphs, in particular, on the Banach space of continuous functions. In Section \ref{sec:contractC} we give a
characterization of generators of Feller semigroups in terms of boundary conditions at vertices of the graph.

The recent reference \cite{Mugnolo:3} is closely related to the results of the present work.
In the main body of the paper we will comment on the related results proved in \cite{Mugnolo:3}.

We mention that the spectrum of generators of semigroups and the spectral mapping are not discussed in the present paper. Also, a class of boundary conditions, the so-called Wentzell boundary conditions, particularly important in the theory of stochastic
processes on metric graphs (see \cite{Freidlin:1}, \cite{KPS2},
\cite{Walsh}), remains out of the scope of the present paper. We
will consider these questions elsewhere.

\subsection*{Acknowledgements} V.~K.~and R.~S.~would like to thank the the Isaac
Newton Institute for Mathematical Sciences for financial support and
hospitality extended to them during their stay in Cambridge in the Spring of
2007. It is also a pleasure to thank the organizers of the programme
``Analysis on graphs and its applications'' for the very inspiring atmosphere there.

\section{Laplace Operators on Metric Graphs. Main Results}\label{sec:basic}

In this section we summarize the terminology used below and present the
main results obtained in the present work.

A finite graph is a 4-tuple $\cG=(V,\cI,\cE,\partial)$, where $V$ is a
finite set of \emph{vertices}, $\cI$ is a finite set of \emph{internal
edges}, $\cE$ is a finite set of \emph{external edges}. Elements in
$\cI\cup\cE$ are called \emph{edges}. The map $\partial$ assigns to each
internal edge $i\in\cI$ an ordered pair of (possibly equal) vertices
$\partial(i):=(v_1,v_2)$ and to each external edge $e\in\cE$ a single
vertex $v$. The vertices $v_1=:\partial^-(i)$ and $v_2=:\partial^+(i)$ are
called the \emph{initial} and \emph{terminal} vertex of the internal edge
$i$, respectively. The vertex $v=\partial(e)$ is the initial vertex of the
external edge $e$. Two vertices $v_1$ and $v_2$ are \emph{adjacent} if there is at least
        one edge $i\in\cI$ with $\partial(i)=(v_1,v_2)$ or
        $\partial(i)=(v_2,v_1)$. If $\partial(i)=(v,v)$, that is,
$\partial^-(i)=\partial^+(i)$ then $i$ is called a \emph{tadpole}. A graph
is called \emph{compact} if $\cE=\emptyset$, otherwise it is
\emph{noncompact}.

Throughout the whole work we will assume that the graph $\cG$ is connected,
that is, for any $v,v^\prime\in V$ there is an ordered sequence $(v_1=v,
v_2,\ldots, v_n=v^\prime)$ such that any two successive vertices in this
sequence are adjacent. In particular, this implies that any vertex of the
graph $\cG$ has nonzero degree, that is, for any vertex there is at least one
edge with which it is incident.

The \emph{degree} $\deg(v)$ of the vertex $v$ is defined as
\begin{equation*}
\deg(v) := |\{e \in \cE\, |\, \partial(e) = v\}| + |\{i\in \cI\,
|\, \partial^-(i) = v\}| + |\{i\in \cI\, |\, \partial^+(i) = v\}|,
\end{equation*}
that is, it is the
number of (internal or external) edges incident with the given
vertex $v$ and by which every tadpole is counted twice.

We will endow the graph with the following metric structure. To each internal
edge $i\in\cI$ we associate an interval $[0,a_i]$ with $a_i>0$
such that the initial vertex of $i$ corresponds to $0$ and the terminal
one to $a_i$. To each external edge $e\in\cE$ we associate the
semiline $[0,+\infty)$. We call the number $a_i$ the length of the internal
edge $i$. We will denote by $\underline{a}$ the vector $(a_i)_{i\in\cI}\in (\R_+)^{|\cI|}$. A
compact or noncompact graph $\cG$ endowed with a metric structure is called
a \emph{metric graph} $(\cG,\underline{a})$.

Given a metric graph $(\cG,\underline{a})$ consider the Hilbert space
\begin{equation}\label{hilbert}
\cH\equiv\cH(\cE,\cI,\underline{a})=\cH_{\cE}\oplus\cH_{\cI},\qquad
\cH_{\cE}=\bigoplus_{e\in\cE}\cH_{e},\qquad
\cH_{\cI}=\bigoplus_{i\in\cI}\cH_{i},
\end{equation}
where $\cH_j=L^2(I_j)$ with
\begin{equation*}
I_j=\begin{cases} [0,a_j] & \text{if}\quad j\in\cI,\\  [0,+\infty) &
\text{if}\quad j\in\cE.\end{cases}
\end{equation*}

By $\cD_j$ with $j\in\cE\cup\cI$ denote the set of all $\psi_j\in\cH_j$
such that $\psi_j$ and its derivative $\psi^\prime_j$ are absolutely
continuous and $\psi^{\prime\prime}_j$ is square integrable. Let
$\cD_j^0$ denote the set of those elements $\psi_j\in\cD_j$ which satisfy
\begin{equation*}
\begin{matrix}
\psi_j(0)=0\\ \psi^\prime_j(0)=0
\end{matrix} \quad \text{for}\quad j\in\cE\qquad\text{and}\qquad
\begin{matrix}
\psi_j(0)=\psi_j(a_j)=0\\
\psi^\prime_j(0)=\psi^\prime_j(a_j)=0
\end{matrix}
\quad\text{for}\quad j\in\cI.
\end{equation*}
Let $\Delta^0$ be the differential operator
\begin{equation}\label{Delta:0}
\left(\Delta^0\psi\right)_j (x) = \frac{d^2}{dx^2} \psi_j(x),\qquad
j\in\cI\cup\cE,\qquad x\in I_j
\end{equation}
with domain
\begin{equation*}
\cD^0=\bigoplus_{j\in\cE\cup\cI} \cD_j^0 \subset\cH.
\end{equation*}
It is straightforward to verify that $\Delta^0$ is a closed symmetric
operator with deficiency indices equal to $|\cE|+2|\cI|$.

Now we begin the discussion of boundary conditions for Laplace operators on metric graphs.
To this end we introduce an auxiliary finite-dimensional Hilbert space
\begin{equation}\label{K:def}
\cK\equiv\cK(\cE,\cI)=\cK_{\cE}\oplus\cK_{\cI}^{(-)}\oplus\cK_{\cI}^{(+)}
\end{equation}
with $\cK_{\cE}\cong\C^{|\cE|}$ and $\cK_{\cI}^{(\pm)}\cong\C^{|\cI|}$. Let
${}^d\cK$ denote the ``double'' of $\cK$, that is, ${}^d\cK=\cK\oplus\cK$.

Let $\cJ\subset \cE\cup\cI$ be a subset of edges, consider $x$ in the
cartesian product $\Bigtimes_{j\in \cJ} I_j$ of these edges, and a function $\psi$
\begin{equation*}
    \psi_{\cJ}(x) = \bigl(\psi_j(x_j),\,j\in \cJ\bigr)^T,
\end{equation*}
where the superscript $T$ denotes transposition. For
\begin{equation*}
    \psi\in\cD:=\bigoplus_{j\in\cE\cup\cI} \cD_j
\end{equation*}
we set
\begin{equation}\label{lin1}
[\psi]:=\underline{\psi}\oplus \underline{\psi}^\prime\in{}^d\cK,
\end{equation}
with $\underline{\psi}$ and $\underline{\psi}^\prime$ defined by
\begin{equation}\label{lin1:add}
\underline{\psi} = \begin{pmatrix}
                        \psi_\cE(0)\\
                        \psi_\cI(0)\\
                        \psi_\cI(\ul{a})\\
                   \end{pmatrix},\qquad
\underline{\psi}' = \begin{pmatrix}
                        \psi_\cE'(0)\\
                        \psi_\cI'(0)\\
                        -\psi_\cI'(\ul{a})\\
                    \end{pmatrix}.
\end{equation}

Let $A$ and $B$ be linear maps of $\cK$ onto itself. By $(A,B)$ we denote
the linear map from ${}^d\cK=\cK\oplus\cK$ to $\cK$ defined by the relation
\begin{equation*}
(A,B)\; (\chi_1\oplus \chi_2) := A\, \chi_1 + B\, \chi_2,
\end{equation*}
where $\chi_1,\chi_2\in\cK$. Set
\begin{equation}\label{M:def}
\cM(A,B) := \Ker\, (A,B).
\end{equation}

The following assumption plays a crucial role throughout the whole work.

\begin{assumption}\label{abcond:neu}
The map $(A,B):\;{}^d\cK\rightarrow\cK$ is surjective, that is, it has
maximal rank equal to $|\cE|+2|\cI|$.
\end{assumption}

Observe that $\Ker A^\dagger \cap \Ker B^\dagger = \{0\}$ under Assumption
\ref{abcond:neu}. Indeed, since the linear map $\begin{pmatrix}A^\dagger \\
B^\dagger\end{pmatrix} = (A,B)^\dagger:\, \cK\rightarrow{}^d\cK$ has maximal
rank equal to $|\cE|+2|\cI|$, it follows that $\Ker \begin{pmatrix}A^\dagger \\
B^\dagger\end{pmatrix}=\{0\}$. Noting that  $\Ker \begin{pmatrix}A^\dagger \\
B^\dagger\end{pmatrix}= \Ker A^\dagger \cap \Ker B^\dagger$ proves the
claim.

\begin{definition}\label{def:equiv}
The boundary conditions $(A,B)$ and $(A',B')$ satisfying Assumption
\ref{abcond:neu} are \emph{equivalent} if the corresponding subspaces $\cM(A,B)$
and $\cM(A',B')$ coincide.
\end{definition}

The boundary conditions $(A,B)$ and $(A',B')$ satisfying Assumption
\ref{abcond:neu} are equivalent if and only if there is an invertible map
$C:\, \cK\rightarrow\cK$ such that $A'= CA$ and $B'=CB$.

Under Assumption \ref{abcond:neu} the inverse $(A+\ii\sk B)^{-1}$ exists
for all $\sk\in\C$ except in a finite subset. Thus,
\begin{equation}\label{uuu:def}
\fS(\sk;A,B):=-(A+\ii\sk B)^{-1} (A-\ii\sk B)
\end{equation}
is well defined for all $\sk\in\C$ but in a finite subset. This
operator plays a central role in the theory of Laplace operators
on metric graphs. In particular,
$\fS(\sk;A,B)=\fS(\sk;A^\prime,B^\prime)$ if and only if
$\cM(A,B)=\cM(A^\prime, B^\prime)$. Hence, we can write
$\fS(\sk;\cM)$ instead of $\fS(\sk;A,B)$ with $\cM=\cM(A,B)$.

With any subspace $\cM\subset{}^d\cK$ of the form \eqref{M:def} we can
associate an extension of $\Delta^{0}$, which is the differential operator
$\Delta(\cM)$ defined by \eqref{Delta:0} with domain
\begin{equation}\label{thru}
\Dom(\Delta(\cM))=\{\psi\in\cD|\; [\psi]\in\cM\}.
\end{equation}
In other words, the domain of the Laplace operator $\Delta(\cM)$ consists
of functions $\psi\in\cD$ satisfying the boundary conditions
\begin{equation}\label{lin2}
 A\underline{\psi}+B\underline{\psi}' = 0,
\end{equation}
with $(A,B)$ subject to \eqref{M:def}. Sometimes we will write
$\Delta(A,B)$ instead of $\Delta(\cM(A,B))$.

Throughout the whole article we adopt the terminology used in
in \cite{Kato} and in \cite{Nagy:Foias}. Recall that the operator
$-\Delta$ is called dissipative if
$$-\Im\langle\psi,\Delta\psi\rangle_{\cH}\geq 0$$ holds for all
$\psi\in\Dom(\Delta)$. The operator $-\Delta$ is accretive if
$$-\Re\langle\psi,\Delta\psi\rangle_{\cH}\geq 0$$ holds for all
$\psi\in\Dom(\Delta)$. A dissipative (respectively accretive)
operator is called maximal, if it does not have a proper
dissipative (respectively accretive) extension. If the domain of a
maximal dissipative (respectively maximal accretive) operator is
dense in $\cH$, then this operator is called m-dissipative
(respectively m-accretive). An m-dissipative or m-accretive
operator is necessarily closed (see \cite{Phillips}, where,
however, a different terminology is used).

Our first main result states that all m-accretive Laplace operators are defined by boundary
conditions satisfying Assumption \ref{abcond:neu}.

\begin{theorem}\label{thm:2:6}
For any m-accretive extension $-\Delta$ of the symmetric operator $-\Delta^0$ the subspace
\begin{equation*}
\cM:=\{[\psi]\, |\, \psi\in\Dom(\Delta) \}\subset{}^d\cK
\end{equation*}
admits the representation $\cM=\Ker(A, B)$ with $(A,B)$ satisfying Assumption \ref{abcond:neu}.
\end{theorem}

By the Lumer-Phillips theorem (see \cite[Theorem II.3.15]{Engel:Nagel} or Theorem IV.4.1 in
\cite{Nagy:Foias}) m-accretive operators are generators of strongly continuous
contraction semigroups, that is, they satisfy
the estimate
\begin{equation*}
\| \e^{t\Delta(A,B)} \| \leq 1,\qquad t>0.
\end{equation*}

Our second main result provides sufficient conditions for the boundary conditions to define an m-accretive operator.

\begin{theorem}\label{thm:accr:neu}
The boundary conditions satisfying Assumption \ref{abcond:neu}
define an m-ac\-cretive Laplace operator $-\Delta(A,B)$ whenever
one of the following equivalent conditions is satisfied
\begin{itemize}
\item[(i)] $\Re(AB^\dagger)\leq 0$;
\item[(ii)] $\fS(\ii\varkappa; A,B)$ defined in \eqref{uuu:def} is a contraction for some (and, thus, for all)
$\varkappa>0$.
\end{itemize}
\end{theorem}

Note that for self-adjoint Laplace operators $-\Delta(A,B)$ this
result has been obtained earlier in \cite{KS9}. We emphasize that
in general the sufficient conditions of Theorem \ref{thm:accr:neu}
need not be necessary. This follows from Example
\ref{ex:4.1} below.

Although Theorem \ref{thm:accr:neu} is stated for differential
operators on graphs, using a concept of the boundary triple (see,
e.g., \cite{Kochubej}), this result can be translated to an
abstract setting, where $\Delta^0$ is replaced by an arbitrary
closed positive symmetric operator on a Hilbert space. A different
description of m-accretive extensions has been obtained by
Tsekanovskii and his coauthors (see \cite{Tsekanovskii},
\cite{Tsekanovskii:2} and references quoted therein).

The methods we use to prove Theorem \ref{thm:accr:neu} can also be
applied to treat m-dissipative extensions of the symmetric
operator $-\Delta^0$. In particular, we obtain a complete
characterization of all m-dissipative extensions, a result which
alternatively can be deduced from Theorem 2 in \cite{Kochubej}.

\begin{theorem}\label{thm:diss:neu}
An extension $-\Delta$ of the symmetric operator $-\Delta^0$ is
m-dissipative if and only if the subspace
\begin{equation*}
\cM:=\{[\psi]\, |\, \psi\in\Dom(\Delta) \}\subset{}^d\cK
\end{equation*}
admits the representation $\cM=\Ker(A, B)$ with $(A,B)$ satisfying
Assumption \ref{abcond:neu} and, in addition, one of the following equivalent
conditions is satisfied
\begin{itemize}
\item[(i)] $\Im(AB^\dagger)\leq 0$;
\item[(ii)] $\fS(-\sk; A,B)$
defined in \eqref{uuu:def} is a contraction for some (and, thus, for all)
$\sk>0$.
\end{itemize}
\end{theorem}

Recall (see \cite{KS1}, \cite{KS3}, and \cite{KS8}) that the extension $-\Delta$ is self-adjoint if and only if the subspace $\cM$ admits the representation $\cM=\Ker(A, B)$ with $(A,B)$ satisfying Assumption \ref{abcond:neu} and, in addition, either
$AB^\dagger$ is self-adjoint or, equivalently, $\fS(-\sk; A,B)$ is unitary.

The proofs of Theorems \ref{thm:2:6}, \ref{thm:accr:neu}, \ref{thm:diss:neu} will be given in Section \ref{sec:diss:accr}.

\subsection{Local boundary Conditions}

With respect to the orthogonal decomposition \eqref{K:def}
any element $\chi$ of
$\cK$ can be represented as a block-vector
\begin{equation}\label{elements}
\chi=\begin{pmatrix}(\chi_e)_{e\in\cE}\\ (\chi^{(-)}_i)_{i\in\cI}\\
(\chi^{(+)}_i)_{i\in\cI}\end{pmatrix}.
\end{equation}
Consider the orthogonal decomposition
\begin{equation}\label{K:ortho}
\cK = \bigoplus_{v\in V} \cL_{v}
\end{equation}
with $\cL_{v}$ the linear subspace of dimension $\deg(v)$ spanned by those
elements \eqref{elements} of $\cK$ which satisfy
\begin{equation}
\label{decomp}
\begin{split}
\chi_e=0 &\quad \text{if}\quad e\in \cE\quad\text{is not incident with the vertex}\quad v,\\
\chi^{(-)}_i=0 &\quad \text{if}\quad v\quad\text{is not an initial vertex of}\quad i\in \cI,\\
\chi^{(+)}_i=0 &\quad \text{if}\quad v\quad\text{is not a terminal vertex
of}\quad i\in \cI.
\end{split}
\end{equation}
Obviously, the subspaces $\cL_{v_1}$ and $\cL_{v_2}$ are orthogonal if
$v_1\neq v_2$.

Set ${^d}\cL_v:=\cL_v\oplus\cL_v\cong\C^{2\deg(v)}$. Obviously, each
$^d\cL_v$ inherits a symplectic structure from ${}^d\cK$ in a canonical way,
such that the orthogonal and symplectic decomposition
\begin{equation}\label{K:decomp}
\bigoplus_{v\in V} {^d}\cL_v = {}^d\cK
\end{equation}
holds.

\begin{definition}\label{propo}
Given the graph $\cG=\cG(V,\cI,\cE,\partial)$, boundary conditions
$(A,B)$ satisfying Assumption \ref{abcond:neu} are called
\emph{local on} $\cG$ if the subspace $\cM(A,B)$ of $\cK$ admits
an orthogonal decomposition
\begin{equation}\label{propo:ortho}
\cM(A,B)=\bigoplus_{v\in V}\;\cM_v,
\end{equation}
where $\cM_v$ are subspaces of ${^d}\cL_v$ of the form \eqref{M:def}
satisfying Assumption \ref{abcond:neu}.

Otherwise the boundary conditions are called \emph{non-local}.
\end{definition}

By Proposition 4.2 in \cite{KS8}, given the graph
$\cG=\cG(V,\cI,\cE,\partial)$, the boundary conditions $(A,B)$ satisfying
Assumption \ref{abcond:neu} are local on $\cG$ if and only if there is an
invertible map $C:\, \cK\rightarrow\cK$ and linear transformations $A_v$
and $B_v$ in $\cL_{v}$ such that the simultaneous orthogonal decompositions
\begin{equation}\label{permut}
CA= \bigoplus_{v\in V} A_v\quad \text{and}\quad CB= \bigoplus_{v\in V} B_v
\end{equation}
are valid such that $\cM_v=\cM(A_v, B_v)$. Alternatively, the boundary conditions $(A,B)$
satisfying Assumption \ref{abcond:neu} are local on $\cG$ if and
only if $\fS(\sk;\cM(A,B))$ admits an orthogonal decomposition
\begin{equation*}
\fS(\sk; \cM(A,B))=\bigoplus_{v\in V}\;\fS(\sk;\cM_v)
\end{equation*}
with respect to \eqref{K:decomp}.

\begin{definition}\label{def:6:2}
A vector $g$ is called \emph{positive} (\emph{respectively strictly positive}),
in symbols $g\succcurlyeq 0$ (respectively $g \succ 0$), if all components of $g$ satisfy $g_j\geq 0$ (respectively $g_j>0$). A $\psi\equiv\{\psi_j\}_{j\in\cI\cup\cE}\in\cH$ is called \emph{positive},
if $\psi(x)\succcurlyeq 0$ for Lebesgue almost all $x$.
A semigroup $\e^{t\Delta(\cM)}$ is called \emph{positivity preserving} if
$\e^{t \Delta(\cM)} \psi$ is positive for all positive
$\psi\in\cH$.
\end{definition}

We say that a $\psi\equiv\{\psi_j\}_{j\in\cI\cup\cE}\in\cH$ is continuous, if $\psi_j(x_j)$ is continuous for all $j\in\cE\cup\cI$ and their boundary values agree at all vertices $v\in V$ with $\deg(v)\geq 2$, that is, for any vertex $v\in V$ with $\deg(v)\geq 2$ there is a number $c_v\in\C$ such that
\begin{equation*}
\psi_j(0) = c_v\qquad\text{for all}\qquad j\in\cE\cup\cI\qquad\text{with}\qquad \partial^-(j)=v
\end{equation*}
and
\begin{equation*}
\psi_j(a_j) = c_v\qquad\text{for all}\qquad j\in\cE\cup\cI\qquad\text{with}\qquad \partial^+(j)=v.
\end{equation*}

\begin{definition}\label{def:cont}
We write
\begin{equation*}
\Dom(\Delta(\cM))\subset \cC(\cG)
\end{equation*}
if all $\psi\in\Dom(\Delta(\cM))$ are continuous. A semigroup $\e^{t\Delta(\cM)}$ is called \emph{continuity preserving} if
$\e^{t \Delta(\cM)} \psi$ is continuous for all continuous
$\psi\in\cH$.
\end{definition}

The following result provides a criterion ensuring that local
boundary conditions define a Laplace operator generating a
strongly continuous contraction semigroup preserving both
continuity and positivity. We set
\begin{equation}\label{h:v:def}
h_v := \begin{pmatrix}1 & 1 & \ldots & 1 \end{pmatrix}^T \in
\cL_v.
\end{equation}
Obviously, $\|h_v\|^2 = \deg(v)$.

\begin{theorem}\label{cor:5:5:neu}
Assume that the graph $\cG$ has no tadpoles. Assume that the
boundary conditions $(A,B)$ are local. The Laplace operator
$-\Delta(A,B)$ generates a strongly continuous contraction
semigroup preserving both continuity and positivity whenever any
of the following equivalent conditions holds:
\begin{itemize}
\item[(i)] {Up to equivalence the boundary conditions $(A_v, B_v)$ are given by
\begin{equation*}
A_v = \1 +\frac{\alpha_v}{\|h_v\|^2} h_v \langle h_v, \cdot\rangle,\qquad B_v = h_v\langle g_v, \cdot\rangle,
\end{equation*}
where $g_v=c h_v$ with $\Re c \leq 0$ if $\alpha_v=0$ and
$c\in\C\setminus\{0\}$ if $\alpha_v=-1$,}
\item[(ii)] {If $\deg(v)\geq 2$, up to equivalence the boundary conditions $(A_v, B_v)$ are given by
\begin{equation*}
A_v=\begin{pmatrix}
    1&-1&0&\ldots&&0&0\\
    0&1&-1&\ldots&&0 &0\\
    0&0&1&\ldots &&0 &0\\
    \vdots&\vdots&\vdots&&&\vdots&\vdots\\
    0&0&0&\ldots&&1&-1\\
    0&0&0&\ldots&&0&-\gamma_v
     \end{pmatrix},\qquad B_v = \begin{pmatrix}
    0&0&0&\ldots&&0&0\\
    0&0&0&\ldots&&0 &0\\
    0&0&0&\ldots &&0 &0\\
    \vdots&\vdots&\vdots&&&\vdots&\vdots\\
    0&0&0&\ldots&&0&0\\
    p&p&p&\ldots&&p&p
     \end{pmatrix},
\end{equation*}
with some $\gamma_v\in \C$, $\Re\gamma_v\geq 0$, $p\geq 0$, and $p\neq
0$ if $\gamma_v=0$.}
\end{itemize}
\end{theorem}

Observe that if $\gamma_v\in \R$, then the boundary
conditions define a self-adjoint Laplace operator
with the so-called $\delta$-type interaction of strength $\gamma_v$
\cite{Exner:Seba}. For $\gamma_v=0$ one has the so-called standard
boundary conditions (see Example 2.6 in \cite{KPS1}).

The proof of Theorem \ref{cor:5:5:neu} is given in Section \ref{sec:pos:preserv}.

The semigroup theory in the Hilbert space $\cH$ can be used to
study the semigroups in other functional spaces on metric graphs.
In particular, our results make it possible to give a complete
characterization of generators of Feller semigroups on graphs with
no internal edges. This result is important in the context of
stochastic processes on metric graphs.

Let $\cC_0(\cG)$ denote the set of all continuous functions on the
graph vanishing at infinity (if $\cE\neq\emptyset$) endowed with
the supremum norm. Obviously, $\cC_0(\cG)$ is a Banach space.
Denote by $\cC_0^2(\cG)$ the subset of $\cC_0(\cG)$ formed by
functions which are twice continuously differentiable on the interior of each edge
of the graph and such that their second derivatives are continuous
at the vertices. Denote by $\CDelta(A,B)$ the differential operator on
$\cC_0(\cG)$ defined by relations similar to \eqref{Delta:0} with
domain
\begin{equation}\label{def:Delta}
\Dom(\CDelta(A,B)) = \{\psi\in\cC_0^2(\cG)\,|\,
A\underline{\psi}+B\underline{\psi}^\prime = 0\}.
\end{equation}
Standard arguments show that $\CDelta(A,B)$ is a closed, densely defined operator.

Following the standard terminology we say that $-\CDelta(A,B)$
generates a Feller semigroup $\e^{t\CDelta(A,B)}$ on $\cC_0(\cG)$ if $\e^{t\CDelta(A,B)}$ is strongly continuous, preserves
positivity, and a contraction with respect to the supremum norm.

\begin{theorem}\label{thm:Feller:star}
Assume that the graph $\cG$ has no internal lines, that is,
$\cI=\emptyset$. Let the boundary conditions $(A,B)$ be local. The
operator $-\CDelta(A,B)$ on $\cC_0(\cG)$ generates a Feller semigroup
if and only if any of the following equivalent conditions holds:
\begin{itemize}
\item[(i)]{Up to equivalence the boundary conditions $(A_v, B_v)$ are given by
\begin{equation}\label{ab:neu}
A_v = \1 +\frac{\alpha_v}{\|h_v\|^2} h_v \langle h_v, \cdot\rangle,\qquad B_v = h_v\langle g_v, \cdot\rangle,
\end{equation}
with some $\alpha_v\in\{0,-1\}$ and some $g_v\in \cL_v$, $g_v\preccurlyeq 0$, subject to the additional restriction $g_v\neq 0$ if $\alpha_v=-1$.}
\item[(ii)] {If $\deg(v)\geq 2$, up to equivalence the boundary conditions $(A_v, B_v)$ are given by
\begin{equation*}
A_v=\begin{pmatrix}
    1&-1&0&\ldots&&0&0\\
    0&1&-1&\ldots&&0 &0\\
    0&0&1&\ldots &&0 &0\\
    \vdots&\vdots&\vdots&&&\vdots&\vdots\\
    0&0&0&\ldots&&1&-1\\
    0&0&0&\ldots&&0&-\gamma_v
     \end{pmatrix},\,\, B_v = \begin{pmatrix}
    0&0&0&\ldots&&0&0\\
    0&0&0&\ldots&&0 &0\\
    0&0&0&\ldots &&0 &0\\
    \vdots&\vdots&\vdots&&&\vdots&\vdots\\
    0&0&0&\ldots&&0&0\\
    p_1&p_2&p_3&\ldots&&p_{n-1}&p_n
     \end{pmatrix},
\end{equation*}
$n=\deg(v)$, with some $\gamma_v\geq 0$,
$p_v=(p_1,p_2,\ldots,p_n)\succcurlyeq 0$ subject to the additional restriction $p_v\neq 0$ if $\gamma_v=0$.}
\end{itemize}
\end{theorem}

The proof of Theorem \ref{thm:Feller:star} is given in Section \ref{sec:contractC}. For general graphs we will prove also a slightly weaker result close
to Theorem \ref{thm:Feller:star} (see Theorem \ref{thm:Feller} below). In \cite{KPS2} we provide a probabilistic proof of Theorem \ref{thm:Feller:star}.

We note that Theorem \ref{thm:Feller:star} is related to a result by Lumer in \cite{Lumer}. Under the assumption that the boundary conditions are given by \eqref{ab:neu} with $\alpha_v=-1$, Theorem 3.1 in \cite{Lumer} states that $-\CDelta(A,B)$ generates a $\cC_0(\cG)$-contraction semigroup if and only if either $g_v\succcurlyeq 0$ or $g_v\preccurlyeq 0$ with $g_v \neq 0$ holds.

\section{Accretive and Dissipative Laplace Operators}\label{sec:diss:accr}

In this section we will prove Theorems \ref{thm:2:6}, \ref{thm:accr:neu}, and \ref{thm:diss:neu}. We start with some auxiliary results.

\begin{lemma}\label{lem:3:1}
The following statements hold under Assumption \ref{abcond:neu}:
\begin{itemize}
\item[(i)]{ $A-\varkappa B$ is invertible for
all $\varkappa > 0$ whenever  $\Re(A B^\dagger) \leq
0$,}
\item[(ii)]{$A-\ii\sk B$ is invertible for all $\sk>0$ whenever $\Im(A B^\dagger) \leq 0$.}
\end{itemize}
\end{lemma}

\begin{proof}
(i) Assume to the contrary that $A-\varkappa B$ is not invertible for some
$\varkappa>0$. Then there is $\chi\in\cK$ such that
\begin{equation}\label{AkB}
(A^\dagger -\varkappa B^\dagger)\chi=0.
\end{equation}
Hence, $AA^\dagger\chi - \varkappa AB^\dagger \chi = 0$. This implies the
equality
\begin{equation*}
\langle\chi, AA^\dagger \chi\rangle  - \varkappa\langle\chi, \Re(A
B^\dagger)\chi\rangle - \ii\varkappa\langle\chi, \Im(AB^\dagger)\chi\rangle = 0.
\end{equation*}
{}From this it follows that
\begin{equation*}
\langle\chi, AA^\dagger \chi\rangle  = \varkappa \langle\chi, \Re(A
B^\dagger)\chi\rangle.
\end{equation*}
Since $\langle\chi, AA^\dagger \chi\rangle \geq 0$ and $\langle\chi, \Re(A
B^\dagger)\chi\rangle \leq 0$, we obtain $\langle\chi, AA^\dagger
\chi\rangle = 0$, which implies $A^\dagger\chi =0$. By \eqref{AkB} we have
$B^\dagger\chi =0$, which contradicts Assumption \ref{abcond:neu}.

(ii) Assume to the contrary that $A-\ii\sk B$ is not invertible for some
$\sk>0$. Then there is $\chi\in\cK$ such that
\begin{equation}\label{AkB:2}
(A^\dagger + \ii\sk B^\dagger)\chi=0.
\end{equation}
Hence, $AA^\dagger\chi + \ii\sk AB^\dagger \chi = 0$. This implies the
equality
\begin{equation*}
\langle\chi, AA^\dagger \chi\rangle  - \sk\langle\chi, \Im(A
B^\dagger)\chi\rangle + \ii\sk\langle\chi, \Re(AB^\dagger)\chi\rangle = 0.
\end{equation*}
{}From this it follows that
\begin{equation*}
\langle\chi, AA^\dagger \chi\rangle  = \sk \langle\chi, \Im(A
B^\dagger)\chi\rangle.
\end{equation*}
Since $\langle\chi, AA^\dagger \chi\rangle \geq 0$ and $\langle\chi, \Im(A
B^\dagger)\chi\rangle \leq 0$, we obtain $\langle\chi, AA^\dagger
\chi\rangle = 0$, which implies $A^\dagger\chi =0$. By \eqref{AkB:2} we have
$B^\dagger\chi =0$, which contradicts Assumption \ref{abcond:neu}.
\end{proof}

\begin{lemma}\label{lem:3:2}
Under Assumption \ref{abcond:neu} the operator $AA^\dagger+BB^\dagger$ is
invertible and the orthogonal projection $P_{\cM^\perp}$ in ${}^d\cK$ onto
the subspace
\begin{equation*}
\cM^\perp :=  \Ran \begin{pmatrix} A^\dagger \\ B^\dagger
\end{pmatrix}
\end{equation*}
orthogonal to the subspace $\cM$ defined in \eqref{M:def}, is given by
\begin{equation}\label{proj:def}
P_{\cM^\perp} = \begin{pmatrix} A^\dagger \\ B^\dagger \end{pmatrix}
(AA^\dagger+BB^\dagger)^{-1} (A, B).
\end{equation}
\end{lemma}

\begin{proof}
Assume that there is a $\chi\in{}^d\cK$ such that
$(AA^\dagger+BB^\dagger)\chi=0$. Then
\begin{equation*}
\langle\chi, (AA^\dagger+BB^\dagger)\chi\rangle = \langle A^\dagger\chi,
A^\dagger\chi\rangle + \langle B^\dagger\chi, B^\dagger\chi\rangle = 0.
\end{equation*}
Thus, $A^\dagger\chi = B^\dagger\chi=0$. Hence, $\chi$ is orthogonal to
both, $\Ran A$ and $\Ran B$, which contradicts Assumption \ref{abcond:neu}.

It is straightforward to verify that \eqref{proj:def} defines an orthogonal
projection. The inclusion $\Ran P_{\cM^\perp}\subset\cM^\perp$ is obvious.
Conversely, a direct calculation shows that
\begin{equation*}
P_{\cM^\perp}\begin{pmatrix}
 A^\dagger \chi\\B^\dagger\chi\end{pmatrix}=\begin{pmatrix}
 A^\dagger \chi\\B^\dagger\chi\end{pmatrix}
\end{equation*}
for any $\chi\in\cK$ and, hence, $\cM^\perp\subset \Ran P_{\cM^\perp}$.
Thus, we have $\Ran P_{\cM^\perp}=\cM^\perp$.
\end{proof}

Assume now that $(A,B)$ satisfies Assumption \ref{abcond:neu} and consider
the Laplace operator $\Delta(\cM)$ corresponding to the subspace $\cM=\cM(A,B)$.
For any $\varphi\in\Dom(\Delta(\cM))$ its quadratic form is given by
\begin{equation*}
\langle\varphi,-\Delta(\cM)\varphi \rangle_{\cH} =
\sum_{j\in\cE\cup\cI}\langle \varphi'_j, \varphi'_j\rangle_{\cH_j}
+ \langle[\varphi],Q[\varphi]\rangle_{{}^d\cK},
\end{equation*}
where $Q=\begin{pmatrix} 0 & \1 \\ 0 & 0 \end{pmatrix}$ with respect to the
orthogonal decomposition ${}^d\cK=\cK\oplus\cK$. Observe that
\begin{equation*}
\begin{split}
\Re \langle\varphi,-\Delta(\cM)\varphi \rangle_{\cH} & = \sum_{j\in\cE\cup\cI}\langle \varphi'_j, \varphi'_j\rangle_{\cH_j} +
\Re\langle[\varphi],Q[\varphi]\rangle_{{}^d\cK},\\
\Im \langle\varphi,-\Delta(\cM)\varphi \rangle_{\cH} & = \Im
\langle[\varphi],Q[\varphi]\rangle_{{}^d\cK}
\end{split}
\end{equation*}
and
\begin{equation*}
\langle[\varphi],Q[\varphi]\rangle_{{}^d\cK} = \langle[\varphi],P_{\cM}QP_{\cM} [\varphi]\rangle_{{}^d\cK},
\end{equation*}
where $P_{\cM}=\1 - P_{\cM^\perp}$.

Thus, we obtain the following result:

\begin{proposition}\label{propo:accret}
Under Assumption \ref{abcond:neu} the operator $-\Delta(\cM)$
\begin{itemize}
\item[(i)] is dissipative
if and only if $\Im P_{\cM} Q P_{\cM} \geq 0$;
\item[(ii)] is accretive whenever
$\Re P_{\cM} Q P_{\cM} \geq 0$.
\end{itemize}
\end{proposition}

The following result establishes a connection between properties of the product $AB^\dagger$ and of the operator $\fS$ defined in \eqref{uuu:def}.

\begin{lemma}\label{lem:accret}
Under Assumption \ref{abcond:neu} the inequality $\Re AB^\dagger\leq 0$
holds if and only if $\|\fS(\ii\varkappa;A,B)\|\leq 1$ for some (and, thus,
for all) $\varkappa>0$. Under the same assumption the inequality $\Im AB^\dagger\leq 0$
holds if and only if $\|\fS(-\sk;A,B)\|\leq 1$ for some (and, thus, for all)
$\sk>0$.
\end{lemma}

\begin{proof}
Assume that $\Re
AB^\dagger\leq 0$.  By Lemma \ref{lem:3:1} $A-\varkappa B$ is invertible for
all $\varkappa>0$, that is, $\fS(\ii\varkappa;A,B)$ is well-defined
by \eqref{uuu:def}. Observe that the boundary conditions $(A,B)$ are
equivalent to the boundary conditions $(A_{\fS}, B_{\fS})$ with
\begin{equation}\label{AfS:BfS:akkr}
A_{\fS} = -\frac{1}{2}(\fS-\1),\qquad B_{\fS} = -\frac{1}{2\varkappa}(\fS+\1),
\end{equation}
where $\fS:=\fS(\ii\varkappa;A,B)$. Indeed, this follows from the equalities
\begin{equation*}
(A-\varkappa B) A_{\fS} = A\qquad\text{and}\qquad (A-\varkappa B) B_{\fS} = B.
\end{equation*}
Therefore, by Sylvester's Inertia Law the inequality $\Re AB^\dagger\leq 0$
holds if and only if $\Re A_{\fS} B_{\fS}^\dagger\leq 0$. Due to
\eqref{AfS:BfS:akkr} we have
\begin{equation*}
A_{\fS} B_{\fS}^\dagger = \frac{1}{4\varkappa}(\fS\fS^\dagger-\1) +
\frac{\ii}{2\varkappa} \Im\fS.
\end{equation*}
Hence, $\Re A_{\fS} B_{\fS}^\dagger\leq 0$ is equivalent to the inequality
$\fS\fS^\dagger\leq \1$. Thus, $\fS^\dagger$ is a contraction and, hence, also its adjoint $\fS$.

Conversely, assume that $\fS(\ii\varkappa_0;A,B)$ is a contraction for some
$\varkappa_0 > 0$. Then the preceding arguments show that $\Re A_{\fS}
B_{\fS}^\dagger \leq 0$, which again by Sylvester's Inertia Law implies $\Re
A B^\dagger \leq 0$. Herewith we also conclude that
$\fS(\ii\varkappa;A,B)$ is a contraction for all $\varkappa > 0$.

We turn to the proof of the second statement. Assume that $\Im AB^\dagger\leq 0$. By Lemma \ref{lem:3:1} $A-\ii\sk B$ is
invertible for all $\sk>0$, that is, $\fS(-\sk;A,B)$ is well-defined by
\eqref{uuu:def}. Observe that the boundary conditions $(A,B)$ are equivalent
to the boundary conditions $(A_{\fS}, B_{\fS})$ with
\begin{equation}\label{AfS:BfS}
A_{\fS} = -\frac{1}{2}(\fS-\1),\qquad B_{\fS} = -\frac{1}{2\ii\sk}(\fS+\1),
\end{equation}
where $\fS:=\fS(-\sk;A,B)$. Indeed, this follows from the equalities
\begin{equation*}
(A-\ii\sk B) A_{\fS} = A\qquad\text{and}\qquad (A-\ii\sk B) B_{\fS} = B.
\end{equation*}
Therefore, by Sylvester's Inertia Law the inequality $\Im AB^\dagger\leq 0$
holds if and only if $\Im A_{\fS} B_{\fS}^\dagger\leq 0$. Due to
\eqref{AfS:BfS} we have
\begin{equation}\label{AfS:BfS:2}
A_{\fS} B_{\fS}^\dagger = \frac{1}{4\ii\sk}(\1-\fS\fS^\dagger) -
\frac{1}{2\sk} \Im\fS.
\end{equation}
Hence, $\Im A_{\fS} B_{\fS}^\dagger\leq 0$ is equivalent to the inequality
$\fS\fS^\dagger\leq \1$. Thus, $\fS(-\sk;A,B)$ is a contraction.

Conversely, assume that $\fS(-\sk_0;A,B)$ is a contraction for some $\sk_0
> 0$. Then the preceding arguments show that $\Im A_{\fS} B_{\fS}^\dagger \leq
0$, which again by Sylvester's Inertia Law implies $\Im A B^\dagger \leq
0$. Herewith we also conclude that $\fS(-\sk;A,B)$ is a contraction for all
$\sk > 0$.
\end{proof}

Using \eqref{proj:def} a simple calculation leads to
\begin{equation*}
P_{\cM^\perp} Q P_{\cM^\perp} = \begin{pmatrix} A^\dagger \\ B^\dagger
\end{pmatrix} (AA^\dagger + BB^\dagger)^{-1} AB^\dagger (AA^\dagger +
BB^\dagger)^{-1}(A, B)
\end{equation*}
such that
\begin{equation}\label{Re:neu}
\Re P_{\cM^\perp} Q P_{\cM^\perp} =\begin{pmatrix} A^\dagger \\ B^\dagger
\end{pmatrix} (AA^\dagger + BB^\dagger)^{-1} \Re(AB^\dagger) (AA^\dagger +
BB^\dagger)^{-1}(A, B)
\end{equation}
and
\begin{equation}\label{Im:neu}
\Im P_{\cM^\perp} Q P_{\cM^\perp} = \begin{pmatrix} A^\dagger \\ B^\dagger
\end{pmatrix} (AA^\dagger + BB^\dagger)^{-1} \Im(AB^\dagger) (AA^\dagger +
BB^\dagger)^{-1}(A, B).
\end{equation}

In the sequel we will need the following lemma with the notation $P^\perp = \1 - P$ for orthogonal projections. We formulate this lemma in the general setting of possibly infinite-dimensional separable Hilbert spaces.

\begin{lemma}\label{lem:geom:neu}
Let $P_1$ and $P_2$ be orthogonal projections in a separable Hilbert space $\fH$.
If the difference $P_1- P_2$ is compact and the pair $(P_1, P_2)$ has vanishing Fredholm index in the sense of \cite{Avron:Seiler},
$\ind(P_1, P_2)=0$,
then the following conditions are equivalent:
\begin{itemize}
\item[(i)]{$P_1 (P_2-P_2^\perp) P_1\geq 0$,}
\item[(ii)]{$P_1^\perp (P_2-P_2^\perp) P_1^\perp\leq 0$.}
\end{itemize}
\end{lemma}

\begin{proof}
\textbf{(i) $\Rightarrow$ (ii).} Observe that
\begin{equation*}
P_1^\perp (P_2-P_2^\perp) P_1^\perp + P_1^\perp = P_1^\perp(P_2-P_1)P_1^\perp - P_1^\perp(P_2^\perp-P_1^\perp)P_1^\perp
\end{equation*}
is compact. Thus, the bounded self-adjoint operator $P_1^\perp (P_2-P_2^\perp) P_1^\perp$ has pure point spectrum.

Assume that (ii) does not hold, that is, the operator $P_1^\perp (P_2-P_2^\perp) P_1^\perp$ has a positive eigenvalue $\lambda>0$.
Denote by $\chi\in \Ran P_1^\perp$ a corresponding eigenvector.

If $P_1(P_2-P_2^\perp)\chi=0$, then we have $(P_2-P_2^\perp)\chi=\lambda\chi$. Hence, $\chi\in\Ran P_2$ and $\lambda=1$. We
arrive at the conclusion $\chi\in \Ran P_1^\perp \cap \Ran P_2$.
Since $\ind(P_1,P_2)=0$, there is a nonzero $\chi^\prime$ lying in $\Ran P_1 \cap \Ran P_2^\perp$. Obviously,
\begin{equation*}
P_1(P_2-P_2^\perp) P_1 \chi^\prime = - \chi^\prime,
\end{equation*}
which contradicts (i).

We turn to the case $P_1(P_2-P_2^\perp)\chi\neq 0$. Then $0\neq\widehat{\chi}\in \Ran P_1$ such that
\begin{equation}\label{make:use}
(P_2-P_2^\perp)\chi = \lambda\chi + \widehat{\chi}.
\end{equation}
This equality implies that
\begin{equation*}
P_1(P_2-P_2^\perp)\chi = \widehat{\chi}.
\end{equation*}
Since $\chi\in\Ran P_1^\perp$, it follows from \eqref{make:use} and $(P_2-P_2^\perp)^2=\1$ that
\begin{equation*}
\begin{split}
P_1(P_2-P_2^\perp)\widehat{\chi} & = P_1(P_2-P_2^\perp)^2 \chi - \lambda P_1(P_2-P_2^\perp)\chi \\
&=P_1\chi - \lambda \widehat{\chi} = - \lambda \widehat{\chi}
\end{split}
\end{equation*}
is valid. Hence, $P_1(P_2-P_2^\perp)P_1\widehat{\chi}=-\lambda\widehat{\chi}$ with $\lambda>0$, which again contradicts (i).

The proof of the implication \textbf{(ii) $\Rightarrow$ (i)} is similar and will, therefore, be omitted.
\end{proof}

\begin{lemma}\label{lem:3:5}
Under Assumption \ref{abcond:neu} the inequality
\begin{equation}\label{3.5b}
\Re P_{\cM} Q P_{\cM} \geq 0
\end{equation}
holds if and only if $\Re(AB^\dagger)\leq 0$. Similarly, the inequality
\begin{equation}\label{3.5a}
\Im P_{\cM} Q P_{\cM} \geq 0
\end{equation}
holds if and only if $\Im(AB^\dagger)\leq 0$.
\end{lemma}

\begin{proof}
\textbf{1.}  We have
\begin{equation}\label{re:Q}
\Re Q = \frac{1}{2}\begin{pmatrix} 0 & \1 \\ \1 & 0 \end{pmatrix} = \frac{1}{2} P_+ - \frac{1}{2} P_-,
\end{equation}
where
\begin{equation}\label{p:strich:pm:def}
P_\pm := \frac{1}{2} \begin{pmatrix} \1 & \mp\1 \\ \mp\1 & \1
\end{pmatrix}
\end{equation}
are orthogonal projections onto the eigenspaces of $\Re Q$, corresponding to
the eigenvalues $\pm\frac{1}{2}$, respectively. It follows from \eqref{Re:neu} and \eqref{re:Q} that
the inequality
\begin{equation}\label{pstrich}
P_{\cM^\perp} (P_+ - P_-) P_{\cM^\perp} \leq 0
\end{equation}
holds if and only if $\Re(AB^\dagger)\leq 0$. Since $P_+$ and $P_\cM$ have equal dimensions, Lemma \ref{lem:geom:neu} can be applied, thus, showing that
inequality \eqref{pstrich} holds if and only if $P_{\cM} (P_+ - P_-) P_{\cM} \geq 0$.

\textbf{2.} We turn to the proof of the second part of the lemma. We have
\begin{equation}\label{im:Q}
\Im Q = \frac{1}{2\ii}\begin{pmatrix} 0 & \1 \\ -\1 & 0 \end{pmatrix} = \frac{1}{2} P'_+ - \frac{1}{2} P'_-,
\end{equation}
where
\begin{equation}\label{p:pm:def}
P'_\pm := \frac{1}{2} \begin{pmatrix} \1 & \mp\ii\1 \\ \pm\ii\1 & \1
\end{pmatrix}
\end{equation}
are orthogonal projections onto the eigenspaces of $\Im Q$, corresponding
to the eigenvalues $\pm\frac{1}{2}$, respectively. {}From \eqref{Im:neu} and \eqref{im:Q} it follows that
the inequality
\begin{equation}\label{pstrich:bis}
P_{\cM^\perp} (P'_+ - P'_-) P_{\cM^\perp} \leq 0
\end{equation}
holds if and only if $\Im(AB^\dagger)\leq 0$. Since $P'_+$ and $P_\cM$ have equal dimensions, by Lemma \ref{lem:geom:neu},
inequality \eqref{pstrich:bis} holds if and only if $P_{\cM} (P'_+ - P'_-) P_{\cM} \geq 0$.

This completes the proof of the lemma.
\end{proof}

For the proof of Theorem \ref{thm:2:6} we need the following lemma.

\begin{lemma}\label{propo:3:6}
An accretive extension $-\Delta(A,B)$ of the symmetric positive
operator $-\Delta^0$ defined in \eqref{Delta:0} is m-accretive if
and only if $(A,B)$ satisfies Assumption \ref{abcond:neu}.
\end{lemma}

\begin{proof}
Assume that the map $(A, B):\, {}^d\cK\rightarrow\cK$ is not
surjective. Then
\begin{equation*}
\Ker\begin{pmatrix} A^\dagger\\ B^\dagger\end{pmatrix} = \Ker
A^\dagger \cap \Ker B^\dagger
\end{equation*}
is nontrivial. A direct calculation shows that the equation
$(-\Delta(A,B)+1)\psi=0$ possesses a solution $\psi\in\cH$ of the
form
\begin{equation}\label{10}
\psi_{j}(x;\sk)=\begin{cases} s_j \e^{- x_j} &\text{for}
                                          \;j\in\cE, \\
                                  \alpha_j \e^{- x_j}+
        \beta_j \e^{ x_j} & \text{for}\; j\in\cI \end{cases}
\end{equation}
if and only if the vectors $s=\{s_e\}_{e\in\cE}\in\cK_{\cE}$,
$\alpha=\{\alpha_i\}_{i\in\cI}\in\cK_{\cI}^{(-)}$, and
$\beta=\{\beta_i\}_{i\in\cI}\in\cK_{\cI}^{(+)}$ satisfy the
homogeneous equation
\begin{equation}\label{11}
Z(A,B)\begin{pmatrix} s\\
                      \alpha\\
                    \beta\end{pmatrix}
                      =0,
\end{equation}
with $Z(A,B):= A X - B Y$, where
\begin{equation}
\label{zet}
X = \begin{pmatrix}\1&0&0\\
                                  0&\1&\1\\
               0&\e^{-\underline{a}}&\e^{+\underline{a}}
               \end{pmatrix}\qquad\text{and}\qquad
Y = \begin{pmatrix}\1&0&0\\
                                  0&\1&-\1\\
               0&-\e^{-\underline{a}}&\e^{+\underline{a}}
               \end{pmatrix}.
\end{equation}
The diagonal $|\cI|\times |\cI|$ matrices $\e^{\pm \underline{a}}$
are given by
\begin{equation}
\label{diag} [\e^{\pm \underline{a}}]_{jk}=\delta_{jk}\e^{\pm
a_{j}}\quad
                       \text{for}\quad j,k\in\;\cI.
\end{equation}
Equation \eqref{11} has indeed a nontrivial solution, since
\begin{equation*}
\Ker Z(A,B)^\dagger = \Ker \left(X^\dagger A^\dagger - Y^\dagger
B^\dagger\right)\supset \Ker A^\dagger \cap \Ker B^\dagger
\end{equation*}
is nontrivial. Thus, $-1$ does not belong to the resolvent set of
$-\Delta(A,B)$ and, hence, $-\Delta(A,B)$ is not m-accretive.

Conversely, assume that $-\Delta(A,B)$ is accretive and $(A,B)$
satisfies Assumption \ref{abcond:neu}. To prove that
$-\Delta(A,B)$ is m-accretive it suffices to show that
$-\Delta(A,B)$ has no proper accretive extensions. Suppose on the
contrary that $-\Delta^\prime$ is a proper accretive extension. Without
loss of generality we can assume that $-\Delta^\prime$ is m-accretive.
Then
\begin{equation*}
\cN:=\{[\psi]\,|\, \psi\in \Dom(-\Delta^\prime)\}\supsetneq
\cM(A,B)=\Ker(A,B).
\end{equation*}
Therefore, there is a pair $(A^\prime,B^\prime)$ such that
$\cN=\Ker(A^\prime,B^\prime)$ and
$\Delta^\prime=\Delta(A^\prime,B^\prime)$. Moreover,
$(A,B):{}^d\cK\rightarrow \cK$ is not surjective. By the preceding
arguments $-\Delta^\prime$ is not m-accretive, a contradiction.
\end{proof}

\begin{proof}[Proof of Theorem \ref{thm:2:6}]
Assume that $-\Delta$ is m-accretive. Consider the linear space
\begin{equation*}
\cM := \left\{[\psi] | \psi\in\Dom(\Delta) \right\} \subset
{}^d\cK.
\end{equation*}
If $\dim \cM \geq |\cE|+ 2|\cI|$, then there exists $A$ and $B$
such that
\begin{equation*}
\cM=\Ker(A,B).
\end{equation*}
Hence, $\Delta=\Delta(A,B)$. By Lemma \ref{propo:3:6}, $(A,B)$
satisfies Assumption \ref{abcond:neu} and, therefore, $\dim \cM =
|\cE|+ 2|\cI|$. Now assume that $\dim \cM < |\cE|+ 2|\cI|$. The
operator $-\Delta$ is m-accretive if and only if its adjoint
$-\Delta^\dagger$ is m-accretive (see \cite[Section
V.3.10]{Kato}). Let us compute the domain of $-\Delta^\dagger$.
For any $\varphi\in\Dom(\Delta)$ and $\psi\in\Dom(\Delta^\dagger)$
we have
\begin{equation*}
\begin{split}
& \langle \psi, -\Delta \varphi\rangle_{\cH}  =
\sum_{j\in\cE\cup\cI} \langle \psi'_j, \varphi'_j\rangle_{\cH_j} +
\langle [\psi], Q [\varphi]\rangle_{{}^d\cK}\\
= & \langle -\Delta^\dagger\psi,  \varphi\rangle_{\cH}  =
\sum_{j\in\cE\cup\cI} \langle \psi'_j, \varphi'_j\rangle_{\cH_j} +
\langle Q [\psi], [\varphi]\rangle_{{}^d\cK}
\end{split}
\end{equation*}
with $Q=\begin{pmatrix} 0 & \1 \\ 0 & 0 \end{pmatrix}$. This
equality implies that
\begin{equation*}
\langle [\psi], \begin{pmatrix} 0 & \1 \\ -\1 & 0 \end{pmatrix}
[\varphi]\rangle_{{}^d\cK} = 0
\end{equation*}
holds for all $\varphi\in\Dom(\Delta)$ and
$\psi\in\Dom(\Delta^\dagger)$. Thus,
\begin{equation*}
\left\{[\psi] | \psi\in\Dom(\Delta^\dagger) \right\} \subset
{}^d\cK
\end{equation*}
is the orthogonal complement in ${}^d\cK$ of the subspace
\begin{equation*}
\left\{\begin{pmatrix} 0 & \1 \\ -\1 & 0 \end{pmatrix} [\varphi]
\Big| \varphi\in\Dom(\Delta)\right\}.
\end{equation*}
Since by the assumption the dimension of this subspace is smaller than
$|\cE|+ 2|\cI|$, we infer
\begin{equation*}
\dim\left\{[\psi] | \psi\in\Dom(\Delta^\dagger) \right\} > |\cE|+
2|\cI|.
\end{equation*}
By the preceding arguments $-\Delta^\dagger$ is not m-accretive,
which is a contradiction.
\end{proof}

\textsc{Proofs of Theorems \ref{thm:accr:neu} and \ref{thm:diss:neu}} are now obtained by combining Proposition \ref{propo:accret}, Lemma \ref{lem:accret}, and Lemma \ref{lem:3:5}.

The following example shows that there are m-accretive Laplace operators $-\Delta(A,B)$ which do not satisfy condition (i) or (ii) in Theorem \ref{thm:accr:neu}.

\begin{figure}[ht]
\centerline{ \unitlength1mm
\begin{picture}(120,40)
\put(20,20){\line(1,0){80}} \put(50,20){\circle*{2}} \put(70,20){\circle*{2}} \put(40,21){$e_1$}
\put(60,21){$i$} \put(80,21){$e_2$} \put(60,20){\vector(1,0){2}}
\end{picture}}
\caption{The graph from Example {\protect\ref{ex:4.1}}. The arrow shows the orientation of the internal edge $i$.}\label{fig:line}
\end{figure}
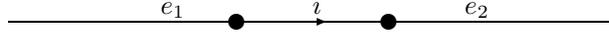

\begin{example}\label{ex:4.1}
On the graph depicted in Fig.\ \ref{fig:line} consider the Laplace operator
$\Delta$ with the boundary conditions
\begin{equation*}
A\underline{\psi}+B \underline{\psi}'=0,
\end{equation*}
where
\begin{equation}\label{AB}
A=\begin{pmatrix} 1 & 0 & -1 & 0 \\ 0 & 1 & 0 & -1 \\ -1 & 0 & 0 & 0 \\ 0 & 1/2 & 0 & 0
\end{pmatrix},\qquad B=\begin{pmatrix} 0 & 0 & 0 & 0 \\ 0 & 0 & 0 & 0 \\
1 & 0 & 1 & 0 \\ 0 & 1 & 0 & 1
\end{pmatrix},
\end{equation}
and where we use the following ordering: $\underline{\psi}=(\psi_{e_1}(0), \psi_{e_2}(0),
\psi_i(0),\psi_i(a))^T$. The boundary conditions \eqref{AB} are, obviously, local in the sense of Definition \ref{propo}. Moreover, $AB^\dagger=B A^\dagger$. Therefore, they define a self-adjoint operator. This operator is unitarily equivalent to the negative of the Laplace operator on the line with two $\delta$-interactions with coupling constants $+1$ and $-1/2$ separated by a distance $a>0$, see \cite[Section II.2.1]{AGHKH}. Hence, the eigenvalues of the operator $-\Delta$ are given by $\lambda=-\varkappa^2$, where $\varkappa>0$ is a solution of the equation
\begin{equation}\label{d:eq}
\left(1+\frac{1}{2\varkappa}\right)\left(-2+\frac{1}{2\varkappa}\right) = \frac{\e^{-2\varkappa a}}{4\varkappa^2}.
\end{equation}
It is easy to verify that for all $a\in (0,1]$ equation \eqref{d:eq} has no positive solutions. Hence, $\spec(-\Delta)=[0,+\infty)$ and so, by the spectral theorem, the operator $-\Delta$ is accretive if $a\in (0,1]$. At the same time the product $A B^\dagger$ does not satisfy the inequality $\Re(A B^\dagger)\leq 0$.
\end{example}

We conclude this section with a description of resolvents of Laplace operators.
The following result is an analogue of Lemma 4.2 in \cite{KS9} for
self-adjoint operators and can be proved in the exactly same way.
The structure of the underlying Hilbert
space $\cH$ \eqref{hilbert} naturally gives rise to the following
definition of integral operators.

Any bounded operator $K$ on the Hilbert space $\cH$ can be represented as a block-operator matrix with entries $K_{j,j^\prime}$ acting from $L^2(I_{j^\prime})$ to $L^2(I_j)$, $j,j^\prime\in\cE\cup\cI$. If all $K_{j,j^\prime}$ are integral operators, we will say that $K$ is an integral operator. More precisely, we adopt the following definition.

\begin{definition}
The operator $K$ on the Hilbert space $\cH$ is called an \emph{integral
operator} if for all $j,j^\prime\in\cE\cup\cI$ there are measurable
functions $K_{j,j^\prime}(\cdot,\cdot)\, :\, I_j\times
I_{j\prime}\rightarrow \C$ with the following properties
\begin{itemize}
\item[(i)]{$K_{j,j^\prime}(x_j,\cdot)\varphi_{j^\prime}(\cdot)\in L^1(I_{j^\prime})$
for almost all $x_j\in I_j$,}
\item[(ii)]{$\psi=K\varphi$ with
\begin{equation}\label{kern}
\psi_j(x_j) = \sum_{j^\prime\in\cE\cup\cI} \int_{I_{j^\prime}}
K_{j,j^\prime}(x_j,y_{j^\prime}) \varphi_{j^\prime}(y_{j^\prime})
dy_{j^\prime}.
\end{equation}}
\end{itemize}
The $(|\cI|+|\cE|)\times(|\cI|+|\cE|)$ matrix-valued function $(x,y):=(x_j, y_{j^\prime})_{j,j^\prime\in \cE\cup\cI}\mapsto
K(x,y)$ with
\begin{equation*}
[K(x,y)]_{j,j^\prime} = K_{j,j^\prime}(x_j,y_{j^\prime})
\end{equation*}
is called the \emph{integral kernel} of the operator $K$.
\end{definition}

Below we will use the following shorthand notation for \eqref{kern}:
\begin{equation*}
\psi(x) = \int^{\cG} K(x,y) \varphi(y) dy.
\end{equation*}

We remark in passing that considering that the integral kernel $K(x,y)$ depends on $\displaystyle x, y\in\Bigtimes_{j\in \cE\cup\cI} I_j$ is consistent with the fact elements of the Hilbert space $\cH$ \eqref{hilbert} are (equivalence classes of) functions
\begin{equation*}
\left(\Bigtimes_{j\in\cE\cup\cI} I_j\right) \ni x \mapsto (\psi_j(x_j))_{j\in\cE\cup\cI} \in \C^{|\cE|+|\cI|}
\end{equation*}
with $\psi_j\in L^2(I_j)$.

There is a different but equivalent way to consider function spaces on metric graphs, which is, in particular, convenient when treating stochastic processes on them \cite{KPS2}. A metric graph can be understood as a disjoint union of intervals $I_j$, where endpoints of $I_j$ and $I_{j^\prime}$ are identified if and only if the edges $j$ and $j^\prime$ are both incident with a vertex $v$. Since this union is a metric space with a natural Lebesgue measure, we may define the space $L^2(\cG)$ as a set of all equivalence classes of square integrable functions $X\mapsto \psi(x)\in\C$, where $x$ belongs to $I_j$ for some $j\in\cE\cup\cI$. Obviously, there is a natural isometric bijection between $\cH$ and $L^2(\cG)$.
In the present article we prefer to work with $\cH$
rather than with $L^2(\cG)$, since this allows for
a presentation of our calculations and results in a convenient and efficient way.

\begin{lemma}\label{lem:Green}
For any subspace $\cM=\cM(A,B)\subset{}^d\cK$ satisfying Assumption
\ref{abcond:neu}, the resolvent
\begin{equation*}
(-\Delta(\cM)-\sk^2)^{-1}\quad\text{for} \quad
\sk^2\in\C\setminus\mathrm{spec}(-\Delta(\cM))\quad\text{with}\quad
\det(A+\ii\sk B)\neq 0,
\end{equation*}
is the integral operator with the $(|\cI|+|\cE|)\times(|\cI|+|\cE|)$
matrix-valued integral kernel $r_{\cM}(x,y;\sk)$, $\Im\sk>0$,
admitting the representation
\begin{equation}\label{r:M:alternativ}
\begin{split}
& r_{\cM}(x,y;\sk)  = r^{(0)}(x,y;\sk)\\ & + \frac{\ii}{2\sk}
\Phi(x,\sk) R_+(\sk;\underline{a})^{-1}[\1-\mathfrak{S}(\sk;\cM)
T(\sk;\underline{a})]^{-1}\mathfrak{S}(\sk;\cM)R_+(\sk;\underline{a})^{-1}
\Phi(y,\sk)^T,
\end{split}
\end{equation}
where $R_+(\sk;\underline{a})$ and $T(\sk;\underline{a})$ are defined by
\begin{equation}\label{U:def:neu}
R_+(\sk;\underline{a}) := \begin{pmatrix} \1 & 0 & 0 \\ 0 & \1 & 0 \\
0 & 0 & \e^{-\ii\sk\underline{a}},\end{pmatrix},\qquad
T(\sk;\underline{a}) := \begin{pmatrix} 0 & 0 & 0 \\ 0 & 0 &
\e^{\ii\sk\underline{a}} \\ 0 & \e^{\ii\sk\underline{a}} & 0
\end{pmatrix}
\end{equation}
with respect to the orthogonal decomposition \eqref{K:def}.
The matrix $\Phi(x,\sk)$ is given by
\begin{equation*}
\Phi(x,\sk) := \begin{pmatrix}\phi(x,\sk) & 0 & 0 \\ 0 & \phi_+(x,\sk) &
\phi_-(x,\sk)
\end{pmatrix}
\end{equation*}
with diagonal matrices $\phi(x,\sk)=\diag\{\e^{\ii\sk x_j}\}_{j\in\cE}$,
$\phi_\pm(x,\sk)=\diag\{\e^{\pm\ii\sk x_j}\}_{j\in\cI}$, $\e^{\pm\ii\sk\underline{a}}=\phi_\pm(\underline{a},\sk)$, and
\begin{equation*}
[r^{(0)}(x,y;\sk)]_{j,j^\prime} = \ii\delta_{j,j^\prime}
\frac{\e^{\ii\sk|x_j-y_j|}}{2\sk},\quad x_j,y_j\in I_j.
\end{equation*}
If $\cI=\emptyset$, representation \eqref{r:M:alternativ} simplifies to
\begin{equation}\label{r:simple}
r_{\cM}(x,y;\sk) = r^{(0)}(x,y;\sk) + \frac{\ii}{2\sk} \phi(x,\sk)
\mathfrak{S}(\sk;\cM) \phi(y,\sk).
\end{equation}
\end{lemma}

The integral kernel $r_{\cM}(x,y;\sk)$ is called
\emph{Green's function} or \emph{Green's matrix}.

\section{Continuity Property}\label{sec:pres:cont}

Let $(A,B)$ be local boundary conditions satisfying Assumption
\ref{abcond:neu}. Obviously, the inclusion
\begin{equation*}
\Dom(\Delta(A,B))\subset \cC(\cG)
\end{equation*}
(see Definition \ref{def:cont}) holds if and only if for any $v\in V^\prime:= \{v\in V| \deg(v)\geq 2\}$ and for any $\begin{pmatrix}\chi_0 \\ \chi_1
\end{pmatrix}\in\Ker (A_v, B_v)$ the vector $\chi_0\in\cK$ is a multiple of
$h_v$, defined in \eqref{h:v:def}.

\begin{theorem}\label{thm:stetig}
For local boundary conditions $(A,B)=\bigoplus_{v\in V} (A_v, B_v)$
satisfying Assumption \ref{abcond:neu} the following statements are
equivalent:
\begin{itemize}
\item[(i)]{$\Dom(\Delta(A,B))\subset \cC(\cG)$;}
\item[(ii)]{For all $v\in V^\prime$ up to equivalence (in the sense of Definition \ref{def:equiv})
the boundary conditions $(A_v,B_v)$ are given by
\begin{equation}\label{Av:Bv}
A_v = \1 + \frac{\alpha_v}{\|h_v\|^2} h_v\langle h_v,\cdot\rangle,\quad
B_v=h_v\langle g_v, \cdot \rangle
\end{equation}
with some $\alpha_v\in\{0,-1\}$ and some $g_v\in\cL_v$ (the case $g_v=0$ is
allowed and corresponds to the Dirichlet boundary conditions) subject to the additional restriction $\langle h_v, g_v\rangle\neq
0$ if $\alpha_v=-1$;}
\item[(iii)]{For all $v\in V^\prime$ up to equivalence (in the sense of Definition \ref{def:equiv})
the boundary conditions $(A_v,B_v)$ are given by
\begin{equation*}
A_v=\begin{pmatrix}
    1&-1&0&\ldots&&0&0\\
    0&1&-1&\ldots&&0 &0\\
    0&0&1&\ldots &&0 &0\\
    \vdots&\vdots&\vdots&&&\vdots&\vdots\\
    0&0&0&\ldots&&1&-1\\
    0&0&0&\ldots&&0&-\gamma_v
     \end{pmatrix},\,\, B_v = \begin{pmatrix}
    0&0&0&\ldots&&0&0\\
    0&0&0&\ldots&&0 &0\\
    0&0&0&\ldots &&0 &0\\
    \vdots&\vdots&\vdots&&&\vdots&\vdots\\
    0&0&0&\ldots&&0&0\\
    p_1&p_2&p_3&\ldots&&p_{n-1}&p_n
     \end{pmatrix},
\end{equation*}
$n=\deg(v)$, with some $\gamma_v\in\C$ and some $p_v=(p_1,
p_2,\ldots,p_n)\in\cL_v$ subject to the additional condition $\langle h_v, p_v\rangle\neq 0$ if $\gamma_v=0$ (the case $p_v=0$, $\gamma_v\neq 0$ is allowed and corresponds to the Dirichlet boundary conditions).}
\end{itemize}
\end{theorem}

\begin{remark}\label{rem:s1:s2}
Equations \eqref{Av:Bv} can be stated equivalently as follows:
\begin{equation}\label{s1:neu}
\mathfrak{S}(\sk; A_v, B_v) = -\1+\frac{2\ii\sk}{1+\ii\sk\langle g_v,
h_v\rangle} h_v\langle g_v, \cdot\rangle\qquad\text{if}\qquad \alpha_v=0
\end{equation}
and
\begin{equation}\label{s2:neu}
\mathfrak{S}(\sk; A_v, B_v) = -\1 + \frac{2}{\langle g_v, h_v\rangle} h_v
\langle g_v,\cdot\rangle,\qquad\text{if}\qquad \alpha_v=-1.
\end{equation}
\end{remark}

Before we turn to the proof of this theorem, we present a simple corollary.

\begin{corollary}\label{cor:stetig:2}
Local boundary conditions $(A,B)=\bigoplus_{v\in V} (A_v, B_v)$
define a self-adjoint Laplace operator with
$\Dom(\Delta(A,B))\subset \cC(\cG)$ if and only if for all $v\in
V^\prime$ up to equivalence (in the sense of Definition
\ref{def:equiv}) the boundary conditions $(A_v,B_v)$ are given by
\begin{equation*}
A_v = \1 + \frac{\alpha_v}{\|h_v\|^2} h_v\langle h_v,\cdot\rangle,\quad B_v= \beta_v h_v\langle h_v, \cdot \rangle
\end{equation*}
with some $\alpha_v\in\{0,-1\}$, $\beta_v\in\R$ and
\begin{equation*}
\Im A_v B_v^\dagger = 0\quad\text{for all}\quad v\in V\setminus V^\prime.
\end{equation*}
The equality $\beta_v = 0$ may only hold if $\alpha_v=0$ (Dirichlet
boundary conditions).

Equivalently, the above statement holds if and only if the
boundary conditions $(A_v,B_v)$ are equivalent either to the
Dirichlet boundary conditions $(\1,0)$ or to the $\delta$-type
boundary conditions with an arbitrary coupling constant $\gamma_v$
as considered in Example 2.6 in \cite{KPS1}.
\end{corollary}

\begin{proof}
The Laplace operator $\Delta(A,B)$ is self-adjoint if and only if $A_v B_v^\dagger$ is self-adjoint for any $v\in V$ (see \cite{KS1}). If $\alpha_v=0$, then $A_v B_v^\dagger$ is
self-adjoint if and only if $B_v$ is self-adjoint. Theorem \ref{thm:stetig}
implies now the claim. If $\alpha_v=-1$, then
\begin{equation*}
A_v B_v^\dagger = g_v\langle h_v,\cdot\rangle - \frac{1}{\|h_v\|^2}\langle h_v, g_v\rangle h_v \langle h_v, \cdot\rangle
\end{equation*}
is self-adjoint if and only if $g_v=\beta_v h_v$ with some real $\beta_v$.
Moreover, $(A_v, B_v)$ has maximal rank if and only if $\beta_v\neq 0$. Again
from Theorem \ref{thm:stetig} the claim follows.
\end{proof}

\begin{remark}
Observe that the boundary conditions referred to in Corollary
\ref{cor:stetig:2} are invariant with respect to permutations of edges. A somewhat related result is Proposition 2.1 in the article \cite{Exner:Turek} by Exner and Turek, which
implies that the $\delta$-type boundary conditions are the only permutation
invariant boundary conditions for which all functions in the domain of
$\Delta(A,B)$ are continuous.
\end{remark}

\medskip

The remainder of this section is devoted to a \textsc{proof of Theorem \ref{thm:stetig}}.

\medskip

\textbf{(i)$\Rightarrow$(ii)}. Consider the subspace ${}^d \cL_v$
associated with an arbitrary vertex $v\in V$. Observe that
$\Dom(\Delta(A_v,B_v))\subset \cC(\cG_v)$ holds if and only if for any
$\begin{pmatrix} \chi_0 \\ \chi_1
\end{pmatrix}\in\Ker (A_v, B_v)$ either $\chi_0=0$ or $\chi_0$ is a nontrivial multiple of
the vector $h_v$. In the second case the equation
\begin{equation}\label{h2}
B_v \chi_1 = - A_v h_v
\end{equation}
has a solution $\chi_1\in\cK$. Hence, we have the following alternative: Either
$B_v=0$ or the subspace
$\cM(A_v, B_v)$ is a linear span of $\{0\}\oplus\Ker B_v$ and $\begin{pmatrix} h_v \\
\chi_1 \end{pmatrix}$, where $\chi_1$ is a solution of \eqref{h2}. If $B_v=0$ we may choose $A_v=\1$ and this corresponds to the Dirichlet boundary conditions at the vertex $v$.

So from now on we will assume that $B_v\neq 0$. Since $\dim\cM(A_v, B_v)=\deg(v)$,
we have $\dim \Ker B_v$ $= \deg(v) - 1$ such that $B_v$ is a rank one operator
and either $\Ker A_v = \{0\}$
or $\dim\Ker A_v = 1$.

First, assume that $\Ker A_v=\{0\}$. Then without loss of generality we can
take $A_v=\1$. Equation \eqref{h2} has a solution if and only if
\begin{equation*}
A_v h_v = h_v \in \Ran B_v.
\end{equation*}
Therefore,
\begin{equation*}
B_v = h_v \langle g_v, \cdot\rangle
\end{equation*}
for some $g_v\in\cL_v$, $g_v\neq 0$.

Second, assume that $\dim\Ker A_v=1$. Then we can take $A_v=P$, an
orthogonal projection of rank $\deg(v)-1$, that is,
\begin{equation*}
P = \1 - f\langle f,\cdot\rangle
\end{equation*}
with some $f\in\cL_v$, $\|f\|=1$. Equation \eqref{h2} has a solution if and
only if $A_v h_v\in\Ran B_v$. Since $B_v$ is of rank one, we have
\begin{equation}\label{b:v:neu}
B_v = \begin{cases}A_v h_v \langle g_v,\cdot\rangle, & \text{if}\quad A_v
h_v\neq 0,\\ \widetilde{g}_v\langle g_v, \cdot\rangle, & \text{if}\quad A_v
h_v = 0\end{cases}
\end{equation}
for some $g_v,\widetilde{g}_v\in\cL_v$, $g_v,\widetilde{g}_v\neq 0$.

If $A_v h_v\neq 0$, then
\begin{equation*}
B_v^\dagger f = g_v\langle A_v h_v, f\rangle = g_v \langle h_v, A_v f\rangle
= 0.
\end{equation*}
Thus, $f\in\Ker A_v^\dagger\cap\Ker B_v^\dagger$. Therefore,
$\dim\Ker\begin{pmatrix} A_v^\dagger \\ B_v^\dagger \end{pmatrix}\geq 1$,
which contradicts Assumption \ref{abcond:neu}. Hence, $A_v h_v = 0$, which
implies that
\begin{equation*}
A_v = \1- \|h_v\|^{-2} h_v\langle h_v,\cdot\rangle\qquad\text{and}\qquad B_v
= \widetilde{g}_v\langle g_v, \cdot\rangle.
\end{equation*}

Assume that $\langle h_v, g_v\rangle =0$. Then $B_v h_v = 0$. Thus, $h_v\in
\Ker A_v\cap \Ker B_v$. This again contradicts Assumption \ref{abcond:neu}.
Thus, $\langle h_v, g_v\rangle \neq 0$.

We claim that the boundary conditions $(A_v, B_v)$ and $(A_v, B_v^\prime)$
with $B_v^\prime = h_v\langle g_v,
\cdot\rangle$ are equivalent, that is, $\cM(A_v, B_v)=\cM(A_v,
B_v^\prime)$. Indeed, let $\chi_0,\chi_1\in\cK$ be an arbitrary solution to
$A_v \chi_0 + B_v\chi_1 = 0$. Then $A_v\chi_0=0$ and $B_v\chi_1=0$.
Therefore, $\langle g_v, \chi_1\rangle = 0$, which implies that $A_v\chi_0 +
B_v^\prime\chi_1=0$. Thus, $\cM(A_v, B_v)\subset\cM(A_v, B_v^\prime)$. Since
$\cM(A_v, B_v)$ and $\cM(A_v, B_v^\prime)$ have equal dimension, we
conclude that $\cM(A_v, B_v)=\cM(A_v, B_v^\prime)$.

\textbf{(ii)$\Rightarrow$(iii)}. Let $(A_v, B_v)$ be given by
\eqref{Av:Bv}. Set
\begin{equation}\label{C:def}
C=\begin{pmatrix}
    1&-1&0&\ldots&&0&0\\
    0&1&-1&\ldots&&0 &0\\
    0&0&1&\ldots &&0 &0\\
    \vdots&\vdots&\vdots&&&\vdots&\vdots\\
    0&0&0&\ldots&&1&-1\\
    1&1&1&\ldots&&1&1
     \end{pmatrix},\qquad
C^\prime=\begin{pmatrix}
    1&-1&0&\ldots&&0&0\\
    0&1&-1&\ldots&&0 &0\\
    0&0&1&\ldots &&0 &0\\
    \vdots&\vdots&\vdots&&&\vdots&\vdots\\
    0&0&0&\ldots&&1&-1\\
    0&0&0&\ldots&&0&1
     \end{pmatrix}.
\end{equation}
A direct calculation shows that $\det C = \deg(v)>0$ and $\det C^\prime = 1$. Obviously,
\begin{equation*}
C h_v = \deg(v) \begin{pmatrix} 0 \\ 0 \\ \vdots \\ 0 \\ 1 \end{pmatrix},\qquad
C^\prime h_v =  \begin{pmatrix} 0 \\ 0 \\ \vdots \\ 0 \\ 1 \end{pmatrix}.
\end{equation*}
If $\alpha_v=-1$, then from \eqref{Av:Bv} it follows that
\begin{equation*}
C A_v = \begin{pmatrix}
    1&-1&0&\ldots&&0&0\\
    0&1&-1&\ldots&&0 &0\\
    0&0&1&\ldots &&0 &0\\
    \vdots&\vdots&\vdots&&&\vdots&\vdots\\
    0&0&0&\ldots&&1&-1\\
    0&0&0&\ldots&&0&0
     \end{pmatrix}
\end{equation*}
and
\begin{equation*}
C B_v = \deg(v) \begin{pmatrix}
    0&0&0&\ldots&&0&0\\
    0&0&0&\ldots&&0 &0\\
    0&0&0&\ldots &&0 &0\\
    \vdots&\vdots&\vdots&&&\vdots&\vdots\\
    0&0&0&\ldots&&0&0\\
    \overline{g_1}&\overline{g_2}&\overline{g_3}&\ldots&&\overline{g_{n-1}}&\overline{g_n}
     \end{pmatrix},
\end{equation*}
where the bar denotes the complex conjugation.
The boundary conditions $(CA_v, CB_v)$ are, obviously, equivalent to those given in (iii) with $\gamma_v=0$ and $p_v=-\deg(v) \overline{g_v}$. If $\alpha_v=0$, then again from \eqref{Av:Bv} it follows that
\begin{equation*}
C^\prime A_v = \begin{pmatrix}
    1&-1&0&\ldots&&0&0\\
    0&1&-1&\ldots&&0 &0\\
    0&0&1&\ldots &&0 &0\\
    \vdots&\vdots&\vdots&&&\vdots&\vdots\\
    0&0&0&\ldots&&1&-1\\
    0&0&0&\ldots&&0&1
     \end{pmatrix}
\end{equation*}
and
\begin{equation*}
C^\prime B_v =  \begin{pmatrix}
    0&0&0&\ldots&&0&0\\
    0&0&0&\ldots&&0 &0\\
    0&0&0&\ldots &&0 &0\\
    \vdots&\vdots&\vdots&&&\vdots&\vdots\\
    0&0&0&\ldots&&0&0\\
    \overline{g_1}&\overline{g_2}&\overline{g_3}&\ldots&&\overline{g_{n-1}}&\overline{g_n}
     \end{pmatrix}.
\end{equation*}
The boundary conditions $(C^\prime A_v, C^\prime B_v)$ are, obviously, equivalent to those given in (iii) with $\gamma_v=-1$ and $p_v=-\overline{g_v}$.

The implication \textbf{(iii)$\Rightarrow$(i)} can be verified by a direct calculation.

This completes the proof of Theorem \ref{thm:stetig}.\hfill $\Box$

\section{Positivity Preserving Contraction Semigroups}\label{sec:pos:preserv}

\subsection{Contraction Semigroups}\label{subsec:contract}

Here we describe those boundary conditions referred to in Theorem \ref{thm:stetig} which
define m-accretive Laplace operators. By the Lumer-Phillips theorem these operators are generators of contraction semigroups.

\begin{lemma}\label{lem:6:1}
Assume that $\Dom(\Delta(A,B))\subset \cC(\cG)$. For local boundary conditions $$(A,B)=\bigoplus_{v\in V} (A_v, B_v)$$
satisfying Assumption \ref{abcond:neu} the following statements are
equivalent for all $v\in V^\prime$:
\begin{itemize}
\item[(i)]{$\mathfrak{S}(\ii\varkappa; A_v, B_v)$ is a contraction for all
$\varkappa>0$,}
\item[(ii)]{the boundary conditions $(A_v,B_v)$ are given by
\begin{equation}\label{Av:Bv:contr:1}
A_v = \1 + \frac{\alpha_v}{\|h_v\|^2} h_v\langle h_v,\cdot\rangle,\quad
B_v=h_v\langle g_v, \cdot \rangle,\qquad \alpha_v\in\{0,-1\},
\end{equation}
where $g_v=c h_v$ with $\Re c \leq 0$ if $\alpha_v=0$ and
$c\in\C\setminus\{0\}$ if $\alpha_v=-1$,}
\item[(iii)]{If $\deg(v)\geq 2$, up to equivalence the boundary conditions $(A_v, B_v)$ are given by
\begin{equation*}
A_v=\begin{pmatrix}
    1&-1&0&\ldots&&0&0\\
    0&1&-1&\ldots&&0 &0\\
    0&0&1&\ldots &&0 &0\\
    \vdots&\vdots&\vdots&&&\vdots&\vdots\\
    0&0&0&\ldots&&1&-1\\
    0&0&0&\ldots&&0&-\gamma_v
     \end{pmatrix},\qquad B_v = \begin{pmatrix}
    0&0&0&\ldots&&0&0\\
    0&0&0&\ldots&&0 &0\\
    0&0&0&\ldots &&0 &0\\
    \vdots&\vdots&\vdots&&&\vdots&\vdots\\
    0&0&0&\ldots&&0&0\\
    p&p&p&\ldots&&p&p
     \end{pmatrix},
\end{equation*}
with some $\gamma_v\in \C$, $\Re\gamma_v\geq 0$, $p\geq 0$, and $p\neq 0$ if $\gamma_v=0$.}
\end{itemize}
\end{lemma}

\begin{proof}
Recall that Theorem \ref{thm:stetig} holds due to the assumption
$\Dom(\Delta(A,B))\subset \cC(\cG)$. Furthermore, $\mathfrak{S}(\ii\varkappa; A_v, B_v)$ is a contraction if and only if
\begin{equation}\label{leq:0}
\mathfrak{S}(\ii\varkappa; A_v, B_v)^\dagger \mathfrak{S}(\ii\varkappa; A_v, B_v)-\1 \leq 0
\end{equation}
in the sense of quadratic forms.

\textbf{(i)$\Leftrightarrow$(ii).} Assume that $\alpha_v=0$. From \eqref{s1:neu} it follows that
\begin{equation}\label{leq:0:bis}
\begin{split}
& \mathfrak{S}(\ii\varkappa; A_v, B_v)^\dagger \mathfrak{S}(\ii\varkappa; A_v, B_v)-\1  =
\frac{4 \varkappa^2 \|h_v\|^2}{|1-\varkappa\langle g_v, h_v\rangle|^2} g_v
\langle g_v, \cdot\rangle\\ &\qquad + \frac{2\varkappa}{1-\varkappa\langle h_v,
g_v\rangle} g_v \langle h_v, \cdot\rangle +
\frac{2\varkappa}{1-\varkappa\langle g_v, h_v\rangle} h_v \langle g_v,
\cdot\rangle.
\end{split}
\end{equation}
Assume that $\mathfrak{S}(\ii\varkappa; A_v, B_v)$ is a contraction. Let $\chi\in\cL_v$ be an arbitrary vector orthogonal to $h_v$. Then, by
\eqref{leq:0},
\begin{equation*}
\begin{split}
& \langle \chi, \left(\mathfrak{S}(\ii\varkappa; A_v, B_v)^\dagger
\mathfrak{S}(\ii\varkappa; A_v, B_v)-\1\right) \chi\rangle  = \frac{4
\varkappa^2 \|h_v\|^2}{|1-\varkappa\langle g_v, h_v\rangle|^2} |\langle g_v,
\chi\rangle|^2 \leq 0
\end{split}
\end{equation*}
and
\begin{equation*}
\begin{split}
& \langle h_v, \left(\mathfrak{S}(\ii\varkappa; A_v, B_v)^\dagger
\mathfrak{S}(\ii\varkappa; A_v, B_v)-\1\right) h_v\rangle  = \frac{4
\varkappa \|h_v\|^2}{|1-\varkappa\langle g_v, h_v\rangle|^2} \Re\langle g_v, h_v\rangle \leq 0.
\end{split}
\end{equation*}
Hence, $g_v = c h_v$ with $\Re c \leq 0$. Conversely, let $g_v = c h_v$ with $\Re c \leq 0$. Then, inequality \eqref{leq:0} follows from \eqref{leq:0:bis}. Thus, $\mathfrak{S}(\ii\varkappa; A_v, B_v)$ is a contraction.

Assume now that $\alpha_v=-1$. From \eqref{s2:neu} it follows that
\begin{equation}\label{leq:0:bis:2}
\begin{split}
\mathfrak{S}(\ii\varkappa; A_v, B_v)^\dagger \mathfrak{S}(\ii\varkappa; A_v,
B_v)-\1 & = \frac{4 \|h_v\|^2}{|\langle g_v, h_v\rangle|^2} g_v \langle g_v,
\cdot\rangle\\ & - \frac{2}{\langle h_v, g_v\rangle} g_v \langle h_v,
\cdot\rangle - \frac{2}{\langle g_v, h_v\rangle} h_v \langle g_v,
\cdot\rangle.
\end{split}
\end{equation}
This implies
\begin{equation*}
\langle h_v, \left(\mathfrak{S}(\ii\varkappa; A_v, B_v)^\dagger
\mathfrak{S}(\ii\varkappa; A_v, B_v)-\1\right) h_v\rangle = 0.
\end{equation*}
If $\mathfrak{S}(\ii\varkappa; A_v, B_v)$ is a contraction, comparing this with \eqref{leq:0}, we conclude that $h_v$ is an eigenvector of
\begin{equation*}
\mathfrak{S}(\sk; A_v, B_v)^\dagger \mathfrak{S}(\sk; A_v, B_v)-\1
\end{equation*}
with eigenvalue zero, that is,
\begin{equation*}
\frac{2\|h_v\|^2}{\langle h_v, g_v\rangle} g_v = 2 h_v,
\end{equation*}
which implies that $g_v$ is a nontrivial multiple of $h_v$, that is, $g_v= c h_v$ with $c\neq 0$. Conversely, if $g_v= c h_v$ with $c\neq 0$, then it follows from \eqref{leq:0:bis:2} that $$\mathfrak{S}(\ii\varkappa; A_v, B_v)^\dagger \mathfrak{S}(\ii\varkappa; A_v,
B_v)=\1.$$ Thus, $\mathfrak{S}(\ii\varkappa; A_v, B_v)$ is a contraction for all $\varkappa>0$.

The equivalence \textbf{(ii)$\Leftrightarrow$(iii)} can be proved in the same way as in Theorem \ref{thm:stetig}.
\end{proof}

Combining Theorem \ref{thm:accr:neu} with Lemma \ref{lem:6:1} we
obtain the following corollary.

\begin{corollary}\label{cor:5:5}
Assume that the boundary conditions $(A,B)$ corresponding to the subspace $\cM=\cM(A,B)$ are local. The Laplace operator $-\Delta(\cM)$ generates a strongly continuous contraction semigroup $\e^{t\Delta(\cM)}$ preserving continuity whenever any of the following equivalent conditions holds:
\begin{itemize}
\item[(i)] {Up to equivalence the boundary conditions $(A_v, B_v)$ are given by
\begin{equation*}
A_v = \1 +\frac{\alpha_v}{\|h_v\|^2} h_v \langle h_v, \cdot\rangle,\qquad B_v = h_v\langle g_v, \cdot\rangle,
\end{equation*}
where $g_v=c h_v$ with $\Re c \leq 0$ if $\alpha_v=0$ and
$c\in\C\setminus\{0\}$ if $\alpha_v=-1$,}
\item[(ii)] {If $\deg(v)\geq 2$, up to equivalence the boundary conditions $(A_v, B_v)$ are given by
\begin{equation*}
A_v=\begin{pmatrix}
    1&-1&0&\ldots&&0&0\\
    0&1&-1&\ldots&&0 &0\\
    0&0&1&\ldots &&0 &0\\
    \vdots&\vdots&\vdots&&&\vdots&\vdots\\
    0&0&0&\ldots&&1&-1\\
    0&0&0&\ldots&&0&-\gamma_v
     \end{pmatrix},\qquad B_v = \begin{pmatrix}
    0&0&0&\ldots&&0&0\\
    0&0&0&\ldots&&0 &0\\
    0&0&0&\ldots &&0 &0\\
    \vdots&\vdots&\vdots&&&\vdots&\vdots\\
    0&0&0&\ldots&&0&0\\
    p&p&p&\ldots&&p&p
     \end{pmatrix},
\end{equation*}
with some $\gamma_v\in \C$, $\Re\gamma_v\geq 0$, $p\geq 0$, and $p\neq 0$ if $\gamma_v=0$.}
\end{itemize}
\end{corollary}

We close this section with an application of our results to evolution equations considered in
\cite{Mugnolo:3}.

\begin{remark}\label{rem:Mugnolo}
We note that Theorem \ref{thm:accr:neu} implies the main part of Corollary 3.3 in \cite{Mugnolo:3}. Indeed, assume that the vertex set $V$ of the graph $\cG$ consists of at least two elements and let $\widetilde{v}\in V$ be arbitrary. Set $\widetilde{V} := V\setminus\{\widetilde{v}\}$. Furthermore, we assume that each vertex $v\in\widetilde{V}$ has degree not smaller than $2$. Consider the following boundary conditions $(A,B)$ on the graph $\cG$. The matrix $A$ is given as a sum $A_1+A_2$. With respect to the orthogonal decomposition \eqref{K:decomp} the matrix $A_1$ is given as a block matrix with blocks
\begin{equation*}
[A_1]_{v,v^\prime} = \begin{cases} E_{v,v^\prime} & \quad\text{if}\quad v,v^\prime\in\widetilde{V},\\
0 & \quad\text{otherwise},  \end{cases}
\end{equation*}
where $E_{v,v^\prime}$ is a $\deg(v)\times\deg(v^\prime)$ matrix of the form
\begin{equation*}
E_{v,v^\prime} = c_{v,v^\prime}\begin{pmatrix}
    0&0&\ldots&&0&0\\
    0&0&\ldots&&0 &0\\
    0&0&\ldots &&0 &0\\
    \vdots&\vdots&&&\vdots&\vdots\\
    0&0&\ldots&&0&0\\
    0&0&\ldots&&0&1
     \end{pmatrix}
\end{equation*}
with $c_{v,v^\prime}\in\C$ arbitrary. The matrices $A_2$ and $B$ are diagonal with respect to this decomposition,
\begin{equation*}
[A_2]_{v,v^\prime} =\begin{cases} \delta_{v,v^\prime}\left(\1-\frac{1}{\|h_v\|^2} h_v\langle h_v, \cdot\rangle\right) & \quad\text{if}\quad v\in\widetilde{V},\\
\delta_{v,v^\prime} \1 & \quad\text{if}\quad v=\widetilde{v}
\end{cases}
\end{equation*}
and
\begin{equation*}
[B]_{v,v^\prime} =\begin{cases} \delta_{v,v^\prime} h_v\langle h_v, \cdot\rangle & \quad\text{if} \quad v\in\widetilde{V},\\
0 & \quad\text{if}\quad v=\widetilde{v}.
\end{cases}
\end{equation*}
It is straightforward to verify that these boundary conditions are equivalent to those given in \cite[Section 2]{Mugnolo:3}. An elementary calculation shows that the inequality $\Re AB^\dagger \leq 0$ holds if and only if the $(|V|-1)\times(|V|-1)$ matrix $C$ with entries $c_{v,v^\prime}$ satisfies $\Re C \geq 0$. Moreover, $AB^\dagger$ is self-adjoint if and only if $C$ is. Thus, from Theorem \ref{thm:accr:neu} combined with the Lumer-Phillips theorem it follows that $-\Delta(A,B)$ generates a contraction semigroup whenever $\Re C \geq 0$ holds. By a result in \cite{KS1} it follows that $-\Delta(A,B)$ generates a self-adjoint semigroup whenever $C$ is self-adjoint.

Obviously, the boundary conditions are local if and only if the matrix $C$ is diagonal.
\end{remark}

\subsection{Positivity Preserving Semigroups}\label{subsec:pos:preserv}

For any matrix $C$ we write $C\succcurlyeq 0$ (respectively, $C\succ
0$) if all entries of the matrix $C$ are nonnegative (respectively,
positive). We write $C_1\succcurlyeq C_2$ (respectively, $C_1\succ C_2$) if
$C_1-C_2\succcurlyeq 0$ (respectively, $C_1-C_2\succ 0$).

\begin{definition}\label{def:positive}
Assume that the boundary conditions $(A,B)$ satisfy Assumption
\ref{abcond:neu}. Set $\cM:=\cM(A,B)$ $=\Ker(A, B)$ according to \eqref{M:def}.
The subspace $\cM\subset{}^d\cK$ is called \emph{positive}, if there is a
$\varkappa_0\geq 0$ such that
$\1+\mathfrak{S}(\ii\varkappa;\cM)\succcurlyeq 0$ for all $\varkappa\geq
\varkappa_0$. It is called \emph{strictly positive}, if
$\1+\mathfrak{S}(\ii\varkappa;\cM)\succ 0$ for all $\varkappa\geq
\varkappa_0$. It is called \emph{locally strictly positive}, if the
boundary conditions defined by $\cM$ are local in the sense of Definition
\ref{propo} and $\1+\mathfrak{S}(\ii\varkappa;\cM_v)\succ 0$ for all
$\varkappa\geq \varkappa_0$ and all $v\in V$. Here $\cM_v$ denotes any
subspace from the orthogonal decomposition \eqref{propo:ortho}.
\end{definition}

In the sequel we will say that boundary conditions $(A,B)$ are positive (respectively strictly positive or locally strictly positive) if the subspace $\cM(A,B)$ is.

We say that a vector $g_v\in\cL_v$ is \emph{sign-definite} if all components
of this vector are either nonnegative or nonpositive, that is, $g_v\succcurlyeq 0$ or $g_v\preccurlyeq 0$. We say that a vector
$g_v\in\cL_v$ is \emph{strictly sign-definite} if it is sign-definite and none
of its components is zero, that is, $g_v\succ 0$ or $g_v\prec 0$.

\begin{proposition}\label{propo:pos}
Assume that the local boundary conditions $(A,B)$ satisfy Assumption
\ref{abcond:neu} and
\begin{equation*}
\Dom(\Delta(A,B))\subset \cC(\cG).
\end{equation*}
These boundary conditions are positive (respectively strictly positive) if and only if for all $v\in V$ up to equivalence the boundary conditions
$(A_v, B_v)$ are given by \eqref{Av:Bv}, where the vector $g_v$ is
sign-definite (respectively strictly
sign-definite).
\end{proposition}

\begin{proof}
If $\deg(v)=1$, the statement is obvious, so let $\deg(v)\geq 2$.

For the case $\alpha_v=0$, by \eqref{s1:neu}, we have
\begin{equation}\label{s1}
\1 + \mathfrak{S}(\ii\varkappa; A_v, B_v) =-\frac{2\varkappa}{1-\varkappa\langle g_v, h_v\rangle} h_v\langle g_v,
\cdot\rangle.
\end{equation}
If $g_v = 0$, then $\1 + \mathfrak{S}(\ii\varkappa; A_v, B_v)\succcurlyeq 0$. Hence, we may assume $g_v\neq 0$. Observe that if $\langle g_v, h_v\rangle=0$, then $g_v$ has strictly
positive as well as strictly negative components such that $h_v\langle g_v,
\cdot\rangle$ has entries of both signs, which implies that
$\1+\fS(\ii\varkappa;\cM)\succcurlyeq 0$ is not valid. Thus, $\langle g_v,
h_v\rangle \neq 0$.

Choose an arbitrary $\varkappa_0 > |\langle g_v, h_v\rangle|^{-1}$. Then
\begin{equation*}
- \sign \langle g_v, h_v\rangle \frac{2\varkappa}{1-\varkappa\langle g_v,
h_v\rangle} > 0
\end{equation*}
for all $\varkappa\geq \varkappa_0$. Hence, the boundary conditions $(A_v,
B_v )$ are positive (strictly positive, respectively) if and only if the
vector $g_v$ is sign-definite (strictly sign-definite, respectively).

Assume now that $\alpha_v=-1$. Recall that by Theorem \ref{thm:stetig}, $\langle g_v, h_v\rangle \neq 0$ and from \eqref{s2:neu} it follows that
\begin{equation}\label{s2}
\1 + \mathfrak{S}(\sk; A_v, B_v) = \frac{2}{\langle g_v, h_v\rangle} h_v
\langle g_v,\cdot\rangle,
\end{equation}
which is independent of $\sk\in\C$. Hence, the boundary conditions $(A_v,
B_v )$ are positive (strictly positive, respectively) if and only if the
vector $g_v$ is sign-definite (strictly sign-definite, respectively).
\end{proof}

By Proposition \ref{propo:pos} the $\delta$-type boundary
conditions (see Example 2.6 in \cite{KPS1}) are locally strictly
positive for all values of the coupling constant $\gamma_v\in\R$.

\begin{corollary}\label{cor:5:3}
Assume that the graph $\cG$ has no tadpoles. Let the boundary
conditions $(A,B)$ be local and satisfy Assumption
\ref{abcond:neu}. Then, $\Dom(\Delta(A,B))\subset \cC(\cG)$ and the
Green's function of $-\Delta(A,B)$ satisfies the inequality
\begin{equation*}
r_{\cM}(x, y; i\varkappa)\succ 0
\end{equation*}
for all sufficiently large $\varkappa>0$,
whenever any of the following equivalent conditions holds:
\begin{itemize}
\item[(i)]{Up to equivalence the boundary conditions $(A_v, B_v)$ are given by
\begin{equation*}
A_v = \1 +\frac{\alpha_v}{\|h_v\|^2} h_v \langle h_v, \cdot\rangle,\qquad B_v = h_v\langle g_v, \cdot\rangle,
\end{equation*}
with some $\alpha_v\in\{0,-1\}$ and some strictly sign-definite
$g_v\in \cL_v$ or $g_v=0$ if $\alpha_v=0$ (Dirichlet boundary
conditions).}
\item[(ii)] {If $\deg(v)\geq 2$, up to equivalence the boundary conditions $(A_v, B_v)$ are given by
\begin{equation*}
A_v=\begin{pmatrix}
    1&-1&0&\ldots&&0&0\\
    0&1&-1&\ldots&&0 &0\\
    0&0&1&\ldots &&0 &0\\
    \vdots&\vdots&\vdots&&&\vdots&\vdots\\
    0&0&0&\ldots&&1&-1\\
    0&0&0&\ldots&&0&-\gamma_v
     \end{pmatrix},\,\, B_v = \begin{pmatrix}
    0&0&0&\ldots&&0&0\\
    0&0&0&\ldots&&0 &0\\
    0&0&0&\ldots &&0 &0\\
    \vdots&\vdots&\vdots&&&\vdots&\vdots\\
    0&0&0&\ldots&&0&0\\
    p_1&p_2&p_3&\ldots&&p_{n-1}&p_n
     \end{pmatrix},
\end{equation*}
$n=\deg(v)$, with some $\gamma_v\geq 0$, $p_v=(p_1,p_2,\ldots,p_n)\succ 0$ or $p_v=0$ if $\gamma_v\neq 0$ (Dirichlet boundary conditions).}
\end{itemize}
\end{corollary}

\begin{proof}
Theorems 5.1 and 6.3 in \cite{KS9} remain valid under the present assumption. Thus, Theorem \ref{thm:stetig} and Proposition \ref{propo:pos} imply the claim.
\end{proof}

Theorem \ref{cor:5:5:neu} follows now from Corollaries \ref{cor:5:5} and \ref{cor:5:3} and
Theorem VI.1.8 in \cite{Engel:Nagel}.

\begin{remark}\label{rem:5:7}
Assume that the graph $\cG$ has no internal lines, that is, $\cI=\emptyset$. Let the boundary
conditions $(A,B)$ be local and satisfy Assumption
\ref{abcond:neu}. Then, $\Dom(\Delta(A,B))\subset \cC(\cG)$ and the
Green's function of $-\Delta(A,B)$ satisfies the inequality
\begin{equation*}
r_{\cM}(x, y; i\varkappa)\succcurlyeq 0
\end{equation*}
for all sufficiently large $\varkappa>0$ if and only if
\begin{equation*}
A_v = \1 +\frac{\alpha_v}{\|h_v\|^2} h_v \langle h_v, \cdot\rangle,\qquad B_v = h_v\langle g_v, \cdot\rangle,
\end{equation*}
with some $\alpha_v\in\{0,-1\}$ and some sign-definite
$g_v\in \cL_v$. To see this, we observe that by \eqref{r:simple} the inequality
\begin{equation*}
[r_{\cM}(x, y; \ii\varkappa)]_{e,e^\prime} \geq 0
\end{equation*}
holds for all $e,e^\prime\in\cE$ with $e\neq e^\prime$ and all $x,y$ if and only if $[\fS(\ii\varkappa;\cM)]_{e,e^\prime}\geq 0$. If $e=e^\prime$ again from \eqref{r:simple} it follows that
\begin{equation}\label{rM:pos}
[r_{\cM}(x, y; \ii\varkappa)]_{e,e} = \frac{1}{2\varkappa}\left(\e^{-\varkappa|x_e-y_e|} +
\e^{-\varkappa(x_e+y_e)} [\fS(\ii\varkappa;\cM)]_{e,e}\right).
\end{equation}
Without loss of generality we can assume that $x_e\geq y_e$. Then the r.h.s.\ of \eqref{rM:pos} can be represented as follows
\begin{equation*}
\frac{\e^{-\varkappa(x_e+y_e)}}{2\varkappa}\left(\e^{2\varkappa y_e} +
 [\fS(\ii\varkappa;\cM)]_{e,e}\right).
\end{equation*}
It is nonnegative for all $y_e\in(0,+\infty)$ if and only if $1+[\fS(\ii\varkappa;\cM)]_{e,e}\geq 0$. Applying Proposition \ref{propo:pos} completes the proof.
\end{remark}


\section{Feller Semigroups}\label{sec:contractC}

In this section we will apply results of the previous sections to
study Feller semigroups on metric graphs and, in particular, we will
prove Theorem \ref{thm:Feller:star}.

Repeating the calculations from the proof of Lemma 4.2 in
\cite{KS9} it is straightforward to verify that
\eqref{r:M:alternativ} is the Green's function of the Laplace operator
$-\CDelta(A,B)$ with domain \eqref{def:Delta}, whenever the boundary conditions satisfy
Assumption \ref{abcond:neu}. Although we will not elaborate on this
observation in detail, the reason for this is the following: The
set $\Dom(\CDelta(A,B))\cap \cH$ is a core for the operator
$\Delta(A,B)$ in the Hilbert space $\cH$ and the closure of
$\CDelta(A,B)$ with respect to the norm of $\cH$ agrees with
$\Delta(A,B)$.

\begin{theorem}\label{thm:Feller}
Assume that the graph $\cG$ has no tadpoles. Let the boundary conditions $(A,B)$ be local. The operator $-\CDelta(A,B)$ on $\cC_0(\cG)$ generates a Feller semigroup whenever any of the following equivalent conditions holds:
\begin{itemize}
\item[(i)]{Up to equivalence the boundary conditions $(A_v, B_v)$ are given by
\begin{equation*}
A_v = \1 +\frac{\alpha_v}{\|h_v\|^2} h_v \langle h_v, \cdot\rangle,\qquad B_v = h_v\langle g_v, \cdot\rangle,
\end{equation*}
with some $\alpha_v\in\{0,-1\}$ and some strictly negative
$g_v\in \cL_v$ or $g_v=0$ if $\alpha_v=0$ (Dirichlet boundary
conditions).}
\item[(ii)] {If $\deg(v)\geq 2$, up to equivalence the boundary conditions $(A_v, B_v)$ are given by
\begin{equation*}
A_v=\begin{pmatrix}
    1&-1&0&\ldots&&0&0\\
    0&1&-1&\ldots&&0 &0\\
    0&0&1&\ldots &&0 &0\\
    \vdots&\vdots&\vdots&&&\vdots&\vdots\\
    0&0&0&\ldots&&1&-1\\
    0&0&0&\ldots&&0&-\gamma_v
     \end{pmatrix},\,\, B_v = \begin{pmatrix}
    0&0&0&\ldots&&0&0\\
    0&0&0&\ldots&&0 &0\\
    0&0&0&\ldots &&0 &0\\
    \vdots&\vdots&\vdots&&&\vdots&\vdots\\
    0&0&0&\ldots&&0&0\\
    p_1&p_2&p_3&\ldots&&p_{n-1}&p_n
     \end{pmatrix},
\end{equation*}
$n=\deg(v)$, with some $\gamma_v\geq 0$, $p_v=(p_1,p_2,\ldots,p_n)\succ 0$ or $p_v=0$ if $\gamma_v\neq 0$ (Dirichlet boundary conditions).}
\end{itemize}
\end{theorem}

\begin{remark}
For the case $\alpha_v=-1$ for all $v\in V$, the integral kernel of $\e^{t\CDelta}$ has been explicitly computed for several graphs in \cite{Okada}.
\end{remark}

\begin{remark}
For boundary conditions considered in Example \ref{rem:Mugnolo} above, Theorem 3.5 in \cite{Mugnolo:3} gives a complete characterization of generators of Feller semigroups.
\end{remark}

For the proof of Theorem \ref{thm:Feller} we need a couple of auxiliary results.

\begin{lemma}\label{lem:6:2}
Assume that the local boundary conditions $(A,B)$
satisfy Assumption \ref{abcond:neu} and the condition $\Dom(\Delta(A,B))\subset
\cC(\cG)$ with $\cM=\cM(A,B)$. Then the following statements are equivalent
\begin{itemize}
\item[(i)]{$\mathfrak{S}(\ii\varkappa; A_v, B_v) h_v \preccurlyeq h_v$ holds for all
$\varkappa>0$,}
\item[(ii)]{Up to equivalence the boundary conditions $(A_v, B_v)$ are given by \eqref{Av:Bv} with $\langle g_v, h_v\rangle \leq 0$
if $\alpha_v=0$ and $\langle g_v, h_v\rangle \neq 0$ if $\alpha_v=-1$,}
\end{itemize}
\end{lemma}

\begin{proof}
Assume that $\alpha_v=0$. Then from \eqref{s1:neu} it follows that
\begin{equation}\label{eq:kitti}
\begin{split}
\mathfrak{S}(\ii\varkappa; A_v, B_v) h_v &= -h_v -
\frac{2\varkappa\langle g_v, h_v\rangle}{1-\varkappa\langle g_v, h_v\rangle}h_v\\
&= h_v \frac{\varkappa\langle g_v, h_v\rangle+1}{\varkappa\langle g_v, h_v\rangle-1}.
\end{split}
\end{equation}
The inequality
\begin{equation}\label{eq:xxx}
\frac{\varkappa\langle g_v, h_v\rangle+1}{\varkappa\langle g_v, h_v\rangle-1} \leq 1
\end{equation}
holds if and only if $\varkappa\langle g_v, h_v\rangle < 1$. Thus, \eqref{eq:xxx} holds for all $\varkappa>0$ if and only if $\langle g_v, h_v\rangle \leq 0$.

Assume that $\alpha_v=-1$. Recall that $\langle g_v, h_v\rangle \neq 0$ by Theorem \ref{thm:stetig}. Then it follows from \eqref{s2:neu} that
\begin{equation*}
\mathfrak{S}(\ii\varkappa; A_v, B_v) h_v = 2 h_v - h_v = h_v.
\end{equation*}
\end{proof}

Note that if $\mathfrak{S}:=\mathfrak{S}(\ii\varkappa; A_v, B_v)\succcurlyeq 0$, the condition (i) in Lemma
\ref{lem:6:2} means that $\mathfrak{S}$ is substochastic. However, under the assumption of this lemma $\mathfrak{S}$ need not be positive. Even the positivity of boundary conditions (cf.~Proposition \ref{propo:pos}) does not imply the positivity of this matrix.

The arguments used in the proof of Lemma \ref{lem:6:2} show also the following result.

\begin{lemma}\label{lem:neu}
Assume that the local boundary conditions $(A,B)$
satisfy Assumption \ref{abcond:neu} and $\Dom(\Delta(\cM))$ $\subset
\cC(\cG)$ with $\cM=\cM(A,B)$. If the inequality
\begin{equation*}
\fS(\ii\varkappa; \cM_v) h_v \preccurlyeq h_v
\end{equation*}
holds for all $\varkappa>0$, then either $\fS(\ii\varkappa; \cM_v) h_v \prec h_v$ or $\fS(\ii\varkappa; \cM_v) h_v = h_v$ holds for all $\varkappa>0$.
\end{lemma}

For the proof it suffices to consider the case $\alpha_v=0$. It follows from Lemma \ref{lem:6:2} that
$\langle g_v, h_v\rangle \leq 0$. Hence
\begin{equation*}
\frac{\varkappa\langle g_v, h_v\rangle+1}{\varkappa\langle g_v, h_v\rangle-1} < 1
\end{equation*}
holds for all $\varkappa>0$. Thus, equation \eqref{eq:kitti} implies $\fS(\ii\varkappa; \cM_v) h_v \prec h_v$.

Without proof we state also the following result, which describes the spectral properties of the matrix $\fS(\ii\varkappa; A_v, B_v)$.

\begin{lemma}\label{lem:lambda}
Assume that the local boundary conditions $(A,B)$
satisfy Assumption \ref{abcond:neu} and $\Dom(\Delta(\cM))$ $\subset
\cC(\cG)$ with $\cM=\cM(A,B)$. Then for all $\varkappa>0$ the spectrum of
$\fS(\ii\varkappa; \cM_v)$ lies in the interval $[-1,1]$.
\end{lemma}

We note that although the spectrum of $\fS(\ii\varkappa; A_v, B_v)$ is real and lies in the interval $[-1,1]$, in general this operator is neither self-adjoint nor a contraction with respect to the $\ell^2$-norm (cf.~Lemma \ref{lem:6:1}).

\begin{lemma}\label{lem:6:invert}
Assume that the local boundary conditions $(A,B)$
satisfy Assumption \ref{abcond:neu} and $\Dom(\Delta(\cM))$ $\subset
\cC(\cG)$ with $\cM=\cM(A,B)$. Then
 $\1-\fS(\ii\varkappa;\cM) T(\ii\varkappa; \underline{a})$ is invertible for all sufficiently large $\varkappa>0$.
\end{lemma}

\begin{proof}
It follows from
\eqref{leq:0:bis}, \eqref{leq:0:bis:2}, Lemmas \ref{lem:6:2} and
\ref{lem:neu} that the norm
$\|\fS(\ii\varkappa;\cM)\|$ is polynomially bounded for all
$\varkappa>0$. Since
\begin{equation*}
\lim_{\varkappa\rightarrow\infty} \|T(\ii\varkappa; \underline{a})\|=0
\end{equation*}
exponentially fast, there is $\varkappa_1\geq 0$ such that $\|\fS(\ii\varkappa;\cM) T(\ii\varkappa; \underline{a})\|<1$ for all $\varkappa>\varkappa_1$.
\end{proof}

Let $h\in\cK$ be the vector with all entries $1$. In particular, we have $h=\bigoplus_{v\in V} h_v$.

\begin{proposition}\label{lem:C:kontr}
Assume that the local boundary conditions $(A,B)$ satisfy
Assumption \ref{abcond:neu}, are strictly positive, and
$\Dom(\Delta(A,B))\subset \cC(\cG)$ with $\cM=\cM(A,B)$. Then the bound
\begin{equation}\label{1:durch:lambda}
\|(-\CDelta(\cM)+\lambda)^{-1}\|_{\infty,\infty} \leq \frac{1}{\lambda}
\end{equation}
holds for all $\lambda>0$ whenever the inequality
\begin{equation}\label{eq:6:5}
\fS(\ii\varkappa; \cM) h \preccurlyeq h
\end{equation}
is valid for all $\varkappa>0$.
\end{proposition}

\begin{proof}
First we observe that it suffices to prove the bound \eqref{1:durch:lambda} for an arbitrary $\lambda>0$. Indeed, assume
that \eqref{1:durch:lambda} holds for some $\lambda_0>0$. Then, by Proposition IV.1.3 in \cite{Engel:Nagel},
\begin{equation}\label{res}
(-\CDelta(\cM)+\lambda)^{-1} = \sum_{n=0}^\infty (\lambda_0-\lambda)^n (-\CDelta(\cM)+\lambda_0)^{-(n+1)}
\end{equation}
holds for all $\lambda\in \C$ satisfying
$|\lambda-\lambda_0|<\|(-\Delta(\cM)+\lambda_0)^{-1}\|^{-1}$. Since, by assumption \eqref{1:durch:lambda},
$\|(-\CDelta(\cM)+\lambda_0)^{-1}\|^{-1}\geq\lambda_0$, we get that the series \eqref{res}
converges for all $\lambda$ satisfying $|\lambda-\lambda_0|<\lambda_0$.
Now, estimating the norm of \eqref{res}, we obtain that
\begin{equation*}
\|(-\CDelta(\cM)+\lambda)^{-1}\|\leq \frac{1}{\lambda_0}\sum_{n=0}^\infty \left(\frac{\lambda_0-\lambda}{\lambda_0}\right)^n = \frac{1}{\lambda}
\end{equation*}
holds for all $\lambda\in(0,2\lambda_0)$. Repeating the above arguments we arrive at the conclusion that \eqref{1:durch:lambda} holds for all $\lambda>0$.

Now we will prove the bound \eqref{1:durch:lambda} for all sufficiently large $\lambda>0$, which
by the preceding argument will imply that \eqref{1:durch:lambda} holds for all $\lambda>0$.

Set $\varkappa=\sqrt{\lambda}$ such that $\sk=\ii\varkappa$ and
$\lambda=\varkappa^2$. Since $r_{\cM}(x,y;\ii\varkappa)\succ 0$ for all sufficiently large $\varkappa>0$,
\begin{equation}\label{eq:norm}
\| (-\CDelta(\cM)+\lambda)^{-1}\|_{\infty,\infty} = \sup_x \left\|
u(x;\varkappa)\right\|_{\ell^\infty(\cE\cup\cI)} = \sup_x
\max_{j\in\cE\cup\cI} u_j(x;\varkappa),
\end{equation}
where
\begin{equation*}
u(x;\varkappa) := \int^{\cG}
r_{\cM}(x,y;\ii\varkappa) 1(y) dy \succcurlyeq 0
\end{equation*}
with $[1(y)]_j = 1$ for all $j\in\cE\cup\cI$.

Consider
\begin{equation}\label{u:0:def}
u^{(0)}(x;\varkappa):=\int^{\cG} r^{(0)}(x,y;\ii\varkappa) 1(y)
dy.
\end{equation}
It is a vector with entries
\begin{equation*}
\frac{1}{2\varkappa} \int_{I_j} \e^{-\varkappa |x_j-y_j|} dy_j, \qquad j\in\cE\cup\cI.
\end{equation*}
An explicit calculation shows that
\begin{equation}\label{j:in:I}
\begin{split}
\frac{1}{2\varkappa} \int_{I_j} \e^{-\varkappa |x_j-y_j|} dy_j &=
 \frac{1}{2\varkappa} \int_0^{x_j} \e^{-\varkappa (x_j-y_j)} dy_j
+
\frac{1}{2\varkappa} \int_{x_j}^{a_j} \e^{-\varkappa (y_j-x_j)} dy_j\\
&=\frac{1}{2\varkappa^2} (1-\e^{-\varkappa x_j}) + \frac{1}{2\varkappa^2}
(1-\e^{-\varkappa (a_j-x_j)})\\
&= \frac{1}{\varkappa^2} - \frac{1}{2\varkappa^2}(\e^{-\varkappa
x_j}+\e^{-\varkappa (a_j-x_j)})
\end{split}
\end{equation}
whenever $j\in\cI$ and
\begin{equation}\label{j:in:E}
\begin{split}
\frac{1}{2\varkappa} \int_{I_j} \e^{-\varkappa |x_j-y_j|} dy_j &
=\frac{1}{2\varkappa} \int_0^{x_j} \e^{-\varkappa (x_j-y_j)} dy_j +
\frac{1}{2\varkappa} \int_{x_j}^{\infty} \e^{-\varkappa (y_j-x_j)} dy_j\\
&=\frac{1}{2\varkappa^2} (1-\e^{-\varkappa x_j}) + \frac{1}{2\varkappa^2}\\
&=\frac{1}{\varkappa^2} -\frac{1}{2\varkappa^2} \e^{-\varkappa
x_j}
\end{split}
\end{equation}
whenever $j\in\cE$.

Now we consider
\begin{equation}\label{u:1:def}
\begin{split}
u^{(1)}(x;\varkappa) & := u(x;\varkappa) - u^{(0)}(x;\varkappa)\\
&\phantom{:}= \int^{\cG}
\left(r(x,y;\ii\varkappa)-r^{(0)}(x,y;\ii\varkappa)\right) 1(y)
dy.
\end{split}
\end{equation}
First we observe that $\chi(\varkappa):=\int^{\cG}
\Phi(y,\ii\varkappa)^T 1(y) dy\in\cK$ is a vector with components
\begin{equation*}
\begin{split}
\int_0^\infty \e^{-\varkappa y_j} dy_j & = \frac{1}{\varkappa},\qquad j\in\cE,\\
\int_0^{a_j} \e^{-\varkappa y_j} dy_j & = \frac{1}{\varkappa} (1-\e^{-\varkappa a_j}),\qquad j\in\cI,\\
\int_0^{a_j} \e^{\varkappa y_j} dy_j & = \frac{1}{\varkappa}
(\e^{\varkappa a_j}-1),\qquad j\in\cI.
\end{split}
\end{equation*}
Therefore,
\begin{equation}\label{ineq:neu}
R_+(\ii\varkappa;\underline{a})^{-1} \chi(\varkappa) = \frac{1}{\varkappa}
(\1-T(\ii\varkappa; \underline{a}))h.
\end{equation}
Hence,
\begin{equation}\label{eq:u1}
\begin{split}
u^{(1)}(x;\varkappa) &=\frac{1}{2\varkappa^2} \Phi(x,\ii\varkappa)
R_+(\ii\varkappa; \underline{a})^{-1} (\1-\fS(\ii\varkappa; \cM)
T(\ii\varkappa; \underline{a}))^{-1} \\ & \qquad\qquad\qquad\cdot\fS(\ii\varkappa; \cM)(\1- T(\ii\varkappa; \underline{a})) h.
\end{split}
\end{equation}

The trivial equality
\begin{equation*}
(\1-\fS(\ii\varkappa; \cM) T(\ii\varkappa; \underline{a})) h
=\fS(\ii\varkappa; \cM) (\1- T(\ii\varkappa; \underline{a})) h
 +
 (\1-\fS(\ii\varkappa; \cM))h
\end{equation*}
and Lemma \ref{lem:6:invert} entail that
\begin{equation}\label{haha}
\begin{split}
h & =(\1-\fS(\ii\varkappa; \cM) T(\ii\varkappa;
\underline{a}))^{-1} \fS(\ii\varkappa; \cM) (\1- T(\ii\varkappa;
\underline{a})) h \\ & + (\1-\fS(\ii\varkappa; \cM)
T(\ii\varkappa; \underline{a}))^{-1} (\1-\fS(\ii\varkappa; \cM))h
\end{split}
\end{equation}
for all sufficiently large $\varkappa>0$. We claim that the inequality
\begin{equation}\label{long:proof}
\Phi(x,\ii\varkappa) R_+(\ii\varkappa;
\underline{a})^{-1}(\1-\fS(\ii\varkappa; \cM) T(\ii\varkappa;
\underline{a}))^{-1} (\1-\fS(\ii\varkappa; \cM))h \succcurlyeq 0
\end{equation}
holds for all large $\varkappa>0$.
Deferring the proof of this inequality to the end of the section we proceed with the proof
of the theorem. Combining \eqref{haha} and \eqref{long:proof} we
arrive at the conclusion
\begin{equation*}
\begin{split}
& \Phi(x,\ii\varkappa) R_+(\ii\varkappa;
\underline{a})^{-1}(\1-\fS(\ii\varkappa; \cM) T(\ii\varkappa;
\underline{a}))^{-1} \fS(\ii\varkappa; \cM) (\1- T(\ii\varkappa;
\underline{a})) h\\ &\qquad \preccurlyeq \Phi(x,\ii\varkappa)
R_+(\ii\varkappa; \underline{a})^{-1} h.
\end{split}
\end{equation*}
Hence, by \eqref{eq:u1}, we obtain the following bounds
\begin{equation*}
u^{(1)}_j(x;\varkappa) \leq \frac{1}{2\varkappa^2}(\e^{-\varkappa
x_j}+\e^{-\varkappa(a_j - x_j)})\qquad\text{if}\qquad j\in\cI
\end{equation*}
and
\begin{equation*}
u^{(1)}_j(x;\varkappa) \leq \frac{1}{2\varkappa^2}\e^{-\varkappa
x_j}\qquad\text{if}\qquad j\in\cE.
\end{equation*}
Combining these bounds with \eqref{j:in:I} and \eqref{j:in:E}, we
see that the inequality $u_j(x;\varkappa)\leq\varkappa^{-2}$
holds for all $j\in\cE\cup\cI$ and all sufficiently large
$\varkappa>0$. Now, from \eqref{eq:norm} the proposition follows.
\end{proof}

Theorem \ref{thm:Feller} now follows immediately from Proposition
\ref{lem:C:kontr} by the Hille-Yosida theorem.

\begin{proof}[Proof of Theorem \ref{thm:Feller:star}]
Due to the Hille-Yosida theorem it suffices to show that the bound
\begin{equation}\label{1:durch:lambda:bis}
\|(-\CDelta(\cM)+\lambda)^{-1}\|_{\infty,\infty} \leq
\frac{1}{\lambda}
\end{equation}
holds for all $\lambda>0$ if and only if the inequality
\begin{equation*}
\fS(\ii\varkappa; \cM) h \preccurlyeq h
\end{equation*}
is valid for all $\varkappa>0$.

To prove this claim we first observe that due to Remark
\ref{rem:5:7} under the present assumptions the Green's function is
positive. Thus, equality \eqref{eq:norm} is valid. The implication
``\eqref{eq:6:5} $\Rightarrow$ \eqref{1:durch:lambda}'' follows from
the arguments used in the proof of Proposition \ref{lem:C:kontr}. (We cannot apply Proposition \ref{lem:C:kontr} directly since the boundary conditions are now assumed to be merely positive rather than strictly positive).
To prove the converse statement we observe that $u^{(0)}$ and $u^{(1)}$ defined in \eqref{u:0:def} and \eqref{u:1:def}, respectively, are given by
\begin{equation*}
u^{(0)}(x;\varkappa) =
\frac{1}{\varkappa^2}h-\frac{1}{2\varkappa^2}\phi(x,\ii\varkappa)
h
\end{equation*}
and
\begin{equation*}
u^{(1)}(x;\varkappa) =
\frac{1}{2\varkappa^2}\phi(x,\ii\varkappa)\fS(\ii\varkappa; \cM) h,
\end{equation*}
where $\phi(x,\ii\varkappa)$ is defined in Lemma \ref{lem:Green}.
Thus,
\begin{equation*}
u_j(x;\varkappa) = \frac{1}{\varkappa^2} + \frac{\e^{-\varkappa
x_j}}{2\varkappa^2} [\fS(\ii\varkappa; \cM) h - h]_j.
\end{equation*}
Now \eqref{1:durch:lambda:bis} with $\lambda=\varkappa^2$ and
\eqref{eq:norm} imply that $\fS(\ii\varkappa; \cM) h -
h\preccurlyeq 0$.
\end{proof}

\subsection{Proof of inequality \textbf{(\ref{long:proof})}}\label{subsec:6:11}

If $\cI=\emptyset$, the proof is trivial and follows directly from
Lemma \ref{lem:6:2}. Thus, we assume further that $\cI\neq\emptyset$. In this case
the proof utilizes the notion of walks on metric graphs (see
\cite{KS8}, \cite{KS9}, \cite{KPS1}). We start with recalling this
notion.

A nontrivial walk $\bw$ on the graph $\cG$ from the edge
$j^\prime\in\cE\cup\cI$ to the edge $j\in\cE\cup\cI$ is a sequence
\begin{equation}\label{walk:def}
(j, v_n, j_n, v_{n-1},\ldots, j_1, v_0, j^\prime)
\end{equation}
such that
\begin{itemize}
\item[(i)]{$j_1,\ldots,j_n\in\cI$;}
\item[(ii)]{the vertices $v_0\in V$ and $v_n\in V$ satisfy $v_0\in\partial(j^\prime)$,
$v_0\in\partial(j_1)$, $v_n\in\partial(j)$, and
$v_n\in\partial(j_n)$;}
\item[(iii)]{for any $k\in\{1,\ldots,n-1\}$ the vertex $v_k\in V$ satisfies $v_k\in\partial(j_k)$ and
$v_k\in\partial(j_{k+1})$;}
\item[(iv)]{$v_k=v_{k+1}$ for some $k\in\{0,\ldots,n-1\}$ if and only if $j_k$ is a tadpole.}
\end{itemize}
If $j,j^\prime\in\cE$ this definition is equivalent to that given in
\cite{KS8}.

The number $n$ appearing in \eqref{walk:def} is the \emph{combinatorial length} $|\bw|_{\mathrm{comb}}$
and the number
\begin{equation*}
|\bw|=\sum_{k=1}^n a_{j_k} >0
\end{equation*}

\noindent is the \emph{metric length} of the walk $\bw$.

A \emph{trivial} walk on the graph $\cG$ from $j^\prime\in\cE\cup\cI$ to
$j\in\cE\cup\cI$ is a triple $(j,v,j^\prime)$ such that
$v\in\partial(j)$ and $v\in\partial(j^\prime)$. In particular, if $\partial(j)=(v_0,v_1)$, then
$(j,v_0,j)$ and $(j,v_1,j)$ are trivial walks, whereas
$(j,v_0,j,v_1,j)$ and $(j,v_1,j,v_0,j)$ are nontrivial walks of
combinatorial length $1$. By convention, both the combinatorial and metric length of a
trivial walk are zero.

We will say that the walk \eqref{walk:def} leaves the edge $j^\prime$ through the
vertex $v_0$ and enters the edge $j$ through the vertex $v_n$. A
trivial walk $(j,v,j^\prime)$ leaves $j^\prime$ and enters $j$ through
the same vertex $v$.

A walk $\bw=(j,v_n,j_n,v_{n-1},\ldots,j_1,v_0,j^\prime)$ \emph{traverses} an
internal edge $i\in\cI$ if $j_k=i$ for some $1\leq k \leq n$. It
\emph{visits} the vertex $v$ if $v_k=v$ for some $0\leq k \leq n$.

We say that the walk \eqref{walk:def} is \emph{transmitted} at the vertex $v_k$, $1\leq k\leq n-1$ if $j_k\neq j_{k+1}$. It is transmitted at the vertex $v_0$ (respectively $v_n$) if $j_1\neq j^\prime$ (respectively $j_n\neq j$). Otherwise the walk is said to be
\emph{reflected}. The walk is called \emph{reflectionless} if
it is transmitted at any vertex visited by this walk.

Under the assumptions of the Theorem \ref{lem:C:kontr}, from
\eqref{leq:0:bis}, \eqref{leq:0:bis:2}, Lemmas \ref{lem:6:2} and
\ref{lem:neu} it follows that the norm
$\|\fS(\ii\varkappa;\cM)\|$ is uniformly bounded for all
$\varkappa>0$. Since
\begin{equation*}
\lim_{\varkappa\rightarrow\infty} \|T(\ii\varkappa;\underline{a})\|=0,
\end{equation*}
there is $\varkappa_1\geq 0$ such that
\begin{equation}\label{eq:series}
(\1 -  \fS(\ii\varkappa;\cM) T(\ii\varkappa;\underline{a}))^{-1} = \sum_{n=0}^\infty
\left(\fS(\ii\varkappa;\cM) T(\ii\varkappa;\underline{a})\right)^n
\end{equation}
converges for all $\varkappa>\varkappa_1$ uniformly in
$\varkappa$. Therefore,
\begin{equation*}
\begin{split}
w(x;\varkappa) &:= \Phi(x,\ii\varkappa) R_+(\ii\varkappa;
\underline{a})^{-1}(\1-\fS(\ii\varkappa; \cM) T(\ii\varkappa;
\underline{a}))^{-1} (\1-\fS(\ii\varkappa; \cM))h\\
&= \Phi(x,\ii\varkappa) R_+(\ii\varkappa;
\underline{a})^{-1} (\1-\fS(\ii\varkappa; \cM))h\\ & + \Phi(x,\ii\varkappa) R_+(\ii\varkappa;
\underline{a})^{-1} \sum_{n=1}^\infty \left(\fS(\ii\varkappa; \cM) T(\ii\varkappa;
\underline{a})\right)^n (\1-\fS(\ii\varkappa; \cM))h.
\end{split}
\end{equation*}
We will now show that $w_j(x;\varkappa)\geq 0$ for all $j\in\cI$. The same statement also holds for $j\in\cE$. Its proof is actually much easier and will, therefore, be omitted.

Let $\cW_{j,j^\prime}^{(\sigma,\sigma^\prime)}$,
$\sigma,\sigma^\prime\in\{+,-\}$ denote the set of all walks from
$j^\prime$ to $j$ leaving the edge $j^\prime$ through the vertex $\partial^{\sigma^\prime}(j^\prime)$ and
entering the edge $j$ through the vertex $\partial^\sigma(j)$.
Observe that for given $j\neq j^\prime$ these four sets are disjoint. For arbitrary
$\sigma\in\{-,+\}$ we will write
\begin{equation*}
\overline{\sigma} := \begin{cases}-, & \text{if}\quad \sigma=+,\\  +, &
\text{if}\quad \sigma=-.\end{cases}
\end{equation*}

We set $V_0:=\{v\in V\,|\, \mathfrak{\fS}(\ii\varkappa;\cM_v)h_v \prec h_v\}$. By Lemma \ref{lem:lambda} the set $V_1:=V\setminus V_0$ agrees with
\begin{equation*}
\{v\in V\,|\, \mathfrak{\fS}(\ii\varkappa;\cM_v)h_v = h_v\}.
\end{equation*}
It is straightforward to verify that
\begin{equation}\label{eq:wj}
\begin{split}
& w_j(x,\varkappa) \\ & \qquad = \e^{-\varkappa x_j} \left[(\1-\fS(\ii\varkappa; \cM))h\right]_{j,-} + \e^{-\varkappa (a_j-x_j)} \left[(\1-\fS(\ii\varkappa; \cM))h\right]_{j,+}\\
&+\sum_{\substack{j^\prime\in\cI\\ \sigma^\prime\in\{+,-\}\\ \partial^{\sigma^\prime}(j^\prime)\in V_0}} \Big(\sum_{\bw\in\cW_{j,j^\prime}^{(-,\overline{\sigma^\prime})}}\e^{-\varkappa x_j} W(\varkappa;\bw) \e^{-\varkappa|\bw|} \e^{-\varkappa a_{j^\prime}} \left[(\1-\fS(\ii\varkappa; \cM))h\right]_{j^\prime,\sigma^\prime}\\
& + \sum_{\bw\in\cW_{j,j^\prime}^{(+,\overline{\sigma^\prime})}}\e^{-\varkappa (a_j-x_j)} W(\varkappa;\bw) \e^{-\varkappa|\bw|} \e^{-\varkappa a_{j^\prime}} \left[(\1-\fS(\ii\varkappa; \cM))h\right]_{j^\prime,\sigma^\prime}\Big)
\end{split}
\end{equation}
holds for all sufficiently large $\varkappa>0$, where the weight $W(\varkappa;\bw)$ associated with the walk $\bw=\{j,v_n,j_n,v_{n-1},\ldots,j_1,v_0,j^\prime\}$ is given by
\begin{equation}\label{W:def}
\begin{split}
W(\varkappa;\bw)  & := [\mathfrak{S}(\ii\varkappa; \cM_v)]_{j,j^\prime}\qquad\text{if}\quad \bw\quad \text{is trivial},\\
W(\varkappa;\bw)  & := [\mathfrak{S}(\ii\varkappa; \cM_v)]_{j,j_1} [\mathfrak{S}(\ii\varkappa,
\cM_{v^\prime})]_{j_1,j^\prime}\qquad\text{if}\quad  n\equiv|\bw|_{\mathrm{comb}} = 1,\\
W(\varkappa;\bw)  & := [\mathfrak{S}(\ii\varkappa; \cM_v)]_{j,j_1}
\displaystyle\left(\prod_{l=1}^{n-1} [\mathfrak{S}(\ii\varkappa;
\cM_{v_l})]_{j_l,j_{l+1}}\right)\\  & \qquad\qquad\cdot [\mathfrak{S}(\ii\varkappa,
\cM_{v^\prime})]_{j_{n},j^\prime}\qquad\text{if}\quad n\equiv|\bw|_{\mathrm{comb}} \geq 2.
\end{split}
\end{equation}

Due to the uniform convergence of the series in \eqref{eq:series}
it suffices to control the leading term in \eqref{eq:wj} only. For
any $j\in\cI$ there are four possible cases
\begin{itemize}
\item[(a)] $\partial^{\pm}(j)\in V_0$,
\item[(b)] $\partial^+(j)\in V_0$ and $\partial^-(j)\in V_1$,
\item[(c)] $\partial^-(j)\in V_0$ and $\partial^+(j)\in V_1$,
\item[(d)] $\partial^{\pm}(j)\in V_1$.
\end{itemize}
We will treat these cases separately.

\emph{Case} (a). By assumption we have $\left[(\1-\fS(\ii\varkappa; \cM))h\right]_{j,\pm}>0$. Therefore
\begin{equation*}
\e^{-\varkappa x_j} \left[(\1-\fS(\ii\varkappa; \cM))h\right]_{j,-} + \e^{-\varkappa (a_j-x_j)} \left[(\1-\fS(\ii\varkappa; \cM))h\right]_{j,+} > 0.
\end{equation*}
Observing that the l.h.s.\ of this inequality is the leading term in \eqref{eq:wj} for large $\varkappa>0$, we arrive at the conclusion that $w_j(x,\varkappa)>0$ for all sufficiently large $\varkappa>0$.

\emph{Case} (b). By assumption we have
\begin{equation*}
\left[(\1-\fS(\ii\varkappa; \cM))h\right]_{j,-}=0\qquad \text{and}\qquad \left[(\1-\fS(\ii\varkappa; \cM))h\right]_{j,+}>0.
\end{equation*}
To determine the leading contribution in \eqref{eq:wj} for large $\varkappa>0$, we introduce the set
\begin{equation}\label{wjminus}
\begin{split}
\cW_j^{(-)} & := \Big\{\bw\in \cW_{j,j^\prime}^{(-,\overline{\sigma^\prime})}\quad\text{for some}\quad j^\prime\in\cI,\,\, j^\prime\neq j\quad\text{and}\quad \sigma^\prime\in\{+,-\}\\ &\quad \text{such that}\quad \partial^{\sigma^\prime}(j^\prime)\in V_0\quad \text{and}\quad |\bw|\leq a_j-a_{j^\prime}\Big\}.
\end{split}
\end{equation}
Obviously, this set may be empty. By $\widetilde{\cW}_j^{(-)}$ we denote the subset of $\cW_j^{(-)}$ formed by the walks with the smallest metric length,
\begin{equation*}
\widetilde{\cW}_j^{(-)} := \left\{ \bw\in \cW_j^{(-)}\quad\text{such that}\quad |\bw|\leq |\bw^\prime|\quad\text{for all}\quad \bw^\prime\in \cW_j^{(-)}\right\}.
\end{equation*}

The leading term in \eqref{eq:wj} is given by
\begin{equation}\label{tri:terms}
\begin{split}
& \e^{-\varkappa (a_j-x_j)} \left[(\1-\fS(\ii\varkappa; \cM))h\right]_{j,+}\\ & + \e^{-\varkappa x_j}\e^{-\varkappa a_j} [\fS(\ii\varkappa; \cM)]_{j,j} \left[(\1-\fS(\ii\varkappa; \cM))h\right]_{j,+}\\ & +
\sum_{\substack{\bw\in \widetilde{\cW}_j^{(-)}\\ \bw\neq\{j,\partial^{-}(j),j\}}}
\e^{-\varkappa x_j} W(\varkappa;\bw) \e^{-\varkappa|\bw|} \e^{-\varkappa a_{j^\prime}} \left[(\1-\fS(\ii\varkappa; \cM))h\right]_{j^\prime,\sigma^\prime},
\end{split}
\end{equation}
where in the last term the walk $\bw$ leaves the edge $j^\prime$ through the vertex $\partial^{\sigma^\prime}(j^\prime)$ and enters the edge $j$ through the vertex $\partial^-(j)$. The second term in this expression corresponds to the trivial walk $\{j,\partial^{-}(j),j\}$. We emphasize that the condition $|\bw|\leq a_j-a_j^\prime$ in \eqref{wjminus} guarantees that the third term is not negligible with respect to the second one.

Observe that the sum of two first terms in \eqref{tri:terms} for all sufficiently large $\varkappa>0$ satisfies the lower bound
\begin{equation*}
\begin{split}
& \e^{-\varkappa (x_j+a_j)} \left[(\1-\fS(\ii\varkappa; \cM))h\right]_{j,+} \left(\e^{2\varkappa x_j}+ [\fS(\ii\varkappa; \cM)]_{j,j}\right)\\
& \geq \e^{-\varkappa (x_j+a_j)} \left[(\1-\fS(\ii\varkappa; \cM))h\right]_{j,+} \left(1+ [\fS(\ii\varkappa; \cM)]_{j,j}\right) > 0
\end{split}
\end{equation*}
uniformly in $x_j\in[0,a_j]$. Here on the last step we used the strict positivity of the boundary conditions (cf.~Definition \ref{def:positive}), which, in particular, implies that $1+ [\fS(\ii\varkappa; \cM_{\partial^+(j)})]_{j,j}$ is strictly positive for all sufficiently large $\varkappa>0$.

We turn to the discussion of the third term in \eqref{tri:terms}.

The following lemma is taken from \cite{KS9}.

\begin{lemma}\label{lem:5:4}
Assume that the graph $\cG$ has no tadpoles. Let $\bw\in
\cW_{j,j^\prime}^{(\sigma,\sigma^\prime)}$ be a walk with the smallest
metric length among all walks in
$\cW_{j,j^\prime}^{(\sigma,\sigma^\prime)}$. Assume that $\bw$ is not
reflectionless. Then there is a reflectionless walk $\bw^\prime\in
\cW_{j,j^\prime}^{(\overline{\sigma},\sigma^\prime)}\cup
\cW_{j,j^\prime}^{(\sigma,\overline{\sigma^\prime})}\cup
\cW_{j,j^\prime}^{(\overline{\sigma},\overline{\sigma^\prime})}$ from $j^\prime$ to
$j$ such that
\begin{equation}\label{relations:alle}
\begin{split}
\mathrm{(i)}\qquad\bw &= \{j,\partial^\sigma(j),\bw^\prime\}\qquad\text{if}\quad
\bw^\prime\in\cW_{j,j^\prime}^{(\overline{\sigma},\sigma^\prime)},\\
\mathrm{(ii)}\qquad\bw &= \{\bw^\prime,\partial^{\sigma^\prime}(j^\prime), j^\prime
\}\qquad\text{if}\quad
\bw^\prime\in\cW_{j,j^\prime}^{(\sigma,\overline{\sigma^\prime})},\\
\mathrm{(iii)}\qquad\bw &=
\{j,\partial^\sigma(j),\bw^\prime,\partial^{\sigma^\prime}(j^\prime),j^\prime\}\qquad\text{if}\quad
\bw^\prime\in\cW_{j,j^\prime}^{(\overline{\sigma},\overline{\sigma^\prime})}.
\end{split}
\end{equation}
\end{lemma}

If all walks $\bw\in \widetilde{\cW}_j^{(-)}$, $\bw\neq(j,\partial^{-}(j),j)$ are reflectionless, then from \eqref{W:def} it follows that the third term in \eqref{tri:terms} is positive. Assume now that a walk $\bw\in \widetilde{\cW}_j^{(-)}$, $\bw\neq(j,\partial^{-}(j),j)$, is not reflectionless. Observe that the possibilities (i) and (iii) in \eqref{relations:alle} cannot occur since in these cases $|\bw|\geq a_j$ which implies that $\bw\notin \cW_j^{(-)}$. In the case (ii) the walk $\bw$ is either of the form
\begin{equation}\label{walk:short}
(j,v_1,j^\prime,v_0,j^\prime)\qquad\text{(if $|\bw|_{\mathrm{comb}}=1$)}
\end{equation}
or
\begin{equation}\label{walk:long}
(j,\ldots,j^{\prime\prime}, v_1,j^\prime,v_0,j^\prime) \qquad\text{(if $|\bw|_{\mathrm{comb}}>1$)}
\end{equation}
with $j^{\prime\prime}\neq j^\prime$.
In both cases $v_1\in V_0$. The walk \eqref{walk:short} enters the edge $j$ through the vertex $\partial^+(j)$ and, hence, does not belong
to $\cW_j^{(-)}$. Thus, $\bw$ is of the form \eqref{walk:long}.
Obviously, the walk
\begin{equation*}
\bw^{\prime\prime} := (j,\ldots,j^{\prime\prime})
\end{equation*}
belongs to $\cW_j^{(-)}$ (since $\bw$ does) and has a metric length strictly smaller than $|\bw|$, which contradicts the assumption that $\bw$ is a shortest walk in $\cW_j^{(-)}$.
Thus, we arrive at the conclusion that the walk $\bw$ is reflectionless. Since the boundary conditions are assumed to be strictly positive, the last term in \eqref{tri:terms} is strictly positive. This proves that $w_j(x,\varkappa)>0$ for all sufficiently large $\varkappa>0$.

\emph{Case} (c) can be handled in the exactly same way.

\emph{Case} (d). Denote by $\widehat{\cW}_j^{(\pm)}$ the set of all walks with the smallest metric length among all walks in the set
\begin{equation}\label{def:set}
\begin{split}
& \Big\{\bw\in \cW_{j,j^\prime}^{(\pm,\overline{\sigma^\prime})}\,\,\text{for some}\,\, j^\prime\in\cI,\, j^\prime\neq j\,\text{and}\,\, \sigma^\prime\in\{+,-\}\\ &\qquad\qquad\qquad
\qquad\qquad\qquad\qquad\qquad\qquad \text{such that}\,\, \partial^{\sigma^\prime}(j^\prime)\in V_0\Big\}.
\end{split}
\end{equation}
By assumption we have
\begin{equation*}
\left[(\1-\fS(\ii\varkappa; \cM))h\right]_{j,\pm}=0.
\end{equation*}
Thus, the leading term in \eqref{eq:wj} is given by
\begin{equation}\label{two:terms}
\begin{split}
& \sum_{\bw\in \widehat{\cW}_j^{(-)}}
\e^{-\varkappa x_j} W(\varkappa;\bw) \e^{-\varkappa|\bw|} \e^{-\varkappa a_{j^\prime}} \left[(\1-\fS(\ii\varkappa; \cM))h\right]_{j^\prime,\sigma^\prime}\\ & +
\sum_{\bw\in \widehat{\cW}_j^{(+)}}
\e^{-\varkappa (a_j-x_j)} W(\varkappa;\bw) \e^{-\varkappa|\bw|} \e^{-\varkappa a_{j^\prime}} \left[(\1-\fS(\ii\varkappa; \cM))h\right]_{j^\prime,\sigma^\prime},
\end{split}
\end{equation}
where $j^\prime$ is the initial edge of the walk $\bw$ and $\partial^{\sigma^\prime}(j^\prime)$ the vertex, through which the walk $\bw$ leaves the edge $j^\prime(\bw)$.

If all walks in $\widehat{\cW}_j^{(-)}$ and $\widehat{\cW}_j^{(+)}$ are reflectionless, then then sum \eqref{two:terms} is positive. Therefore, we assume that there is a non-reflectionless walk $\bw$ belonging, say, to $\widehat{\cW}_j^{(-)}$. Observe that (ii) and (iii) in \eqref{relations:alle} cannot occur. Indeed, this would contradict to the assumption that $\bw$ is a walk with the smallest metric length among all walks in the set \eqref{def:set}.

Thus, by Lemma \ref{lem:5:4}, there is a reflectionless walk $\bw^\prime\in \cW_{j,j^\prime}^{(+,\overline{\sigma^\prime})}$ such that
\begin{equation}\label{www}
\bw = (j,\partial^-(j),\bw^\prime).
\end{equation}
We claim that $\bw^\prime\in\widehat{\cW}_j^{(+)}$. Assume to the contrary that there is a walk $\bw^{\prime\prime}\in\widehat{\cW}_j^{(+)}$ with $|\bw^{\prime\prime}| < |\bw^\prime|$. Then, the walk $\{j,\partial^-(j),\bw^{\prime\prime}\}$ has a metric length
\begin{equation*}
|\bw^{\prime\prime}| + a_j < |\bw^\prime| + a_j = |\bw|.
\end{equation*}
Thus, the walk $\{j,\partial^-(j),\bw^{\prime\prime}\}$ belongs to $\widehat{\cW}_j^{(-)}$ and has a length smaller than $|\bw|$. Since $\bw\in\widehat{\cW}_j^{(-)}$, this is a contradiction.

Observe that from \eqref{W:def} and \eqref{www} it follows that
\begin{equation*}
W(\varkappa;\bw) = [\fS(\ii\varkappa; \cM_{\partial^-(j)})]_{j,j} W(\varkappa;\bw^\prime)\qquad\text{and}\qquad |\bw|=|\bw^\prime| + a_j,
\end{equation*}
where $W(\varkappa;\bw^\prime)>0$.
Hence, the sum of the contributions of the walks $\bw$ and $\bw^\prime$ to \eqref{two:terms} is given by
\begin{equation}\label{rhs}
\begin{split}
&
\e^{-\varkappa x_j} W(\varkappa;\bw) \e^{-\varkappa|\bw|} \e^{-\varkappa a_{j^\prime}} \left[(\1-\fS(\ii\varkappa; \cM))h\right]_{j^\prime,\sigma^\prime}\\ & +
\e^{-\varkappa (a_j-x_j)} W(\varkappa;\bw^\prime) \e^{-\varkappa|\bw^\prime|} \e^{-\varkappa a_{j^\prime}} \left[(\1-\fS(\ii\varkappa; \cM))h\right]_{j^\prime,\sigma^\prime}\\=\,\,& \e^{-\varkappa x_j} \e^{-\varkappa (a_j+a_{j^\prime})} \e^{-\varkappa|\bw^\prime|} \left[(\1-\fS(\ii\varkappa; \cM))h\right]_{j^\prime,\sigma^\prime}\\ & \qquad\qquad \cdot \left(\e^{2\varkappa x_j}+[\fS(\ii\varkappa; \cM_{\partial^-(j)})]_{j,j}\right) W(\varkappa;\bw^\prime)\\ \geq \,\,& \e^{-\varkappa x_j} \e^{-\varkappa (a_j+a_{j^\prime})} \e^{-\varkappa|\bw^\prime|} \left[(\1-\fS(\ii\varkappa; \cM))h\right]_{j^\prime,\sigma^\prime}\\ & \qquad\qquad \cdot \left(1+[\fS(\ii\varkappa; \cM_{\partial^-(j)})]_{j,j}\right) W(\varkappa;\bw^\prime)
\end{split}
\end{equation}
uniformly in $x_j\in [0,a_j]$. The strict positivity of the boundary conditions (cf.~Definition \ref{def:positive}) implies that $1+ [\fS(\ii\varkappa; \cM_{\partial^-(j)})]_{j,j}$ is strictly positive for all sufficiently large $\varkappa>0$.
Thus, the r.h.s.\ of \eqref{rhs} is strictly positive for all sufficiently large $\varkappa>0$.

For any other non-reflectionless $\widetilde{\bw}\in \widehat{\cW}_j^{(-)}$ there is a reflectionless walk $\widetilde{\bw}^\prime\in\cW_{j,j^\prime}^{(+,\overline{\sigma^\prime})}$ such that $\widetilde{\bw} = (j,\partial^-(j),\widetilde{\bw}^\prime)$. Obviously, $\bw^\prime$ and $\widetilde{\bw}^\prime$ are different. Thus, $w_j(x,\varkappa)>0$ for all sufficiently large $\varkappa>0$.

This completes the proof of inequality \eqref{long:proof}.\hfill $\Box$



\begin{thebibliography}{50}

\bibitem{AGHKH} S.~Albeverio, F.~Gesztesy, R.~H{\o}egh-Krohn, and H.~Holden,
\textit{Solvable Models in Quantum Mechanics}, Springer, Berlin, 1988.

\bibitem{Ali:Mehmeti} F.~Ali Mehmeti, \textit{Nonlinear Waves in Networks}, Mathematical Research Vol.~80. Berlin, Akademie Verlag, 1994.

\bibitem{Angad-Guar} H.~W.~K.~Angad-Gaur, B.~Gaveau, and M.~Okada, \textit{Explicit heat
kernel on generalized cones}, SIAM J. Math. Anal. \textbf{25} (1994), 1562
-- 1576.

\bibitem{Avron:Seiler} J.~Avron, R.~Seiler, and B.~Simon, \textit{The index of a pair of projections}, J. Funct. Anal.
\textbf{120} (1994), 220 -- 237.

\bibitem{Badalin:1} A.~Badanin, J.~Br\"{u}ning, E.~Korotyaev, and I.~Lobanov \textit{Schr\"{o}dinger operators on armchair nanotubes. I}, Preprint \texttt{arXiv:0707.3909} (2007).

\bibitem{Badalin:2} A.~Badanin, J.~Br\"{u}ning, and E.~Korotyaev, \textit{Schr\"{o}dinger operators on armchair nanotubes. II}, Preprint \texttt{arXiv:0707.3900} (2007).

\bibitem{Below} J.~von Below, \textit{A maximum principle for semilinear parabolic network equations}, in J.~A.~Goldstein, F.~Kappel, and W.~Schappacher (eds.),
\textit{Differential equations with applications in biology, physics, and engineering}, (Proceedings of the international conference, held in Leibnitz, Austria),
Lecture Notes in Pure and Applied Mathematics Vol.~133. New York , Marcel Dekker, 1991. pp.~37 -- 45.

\bibitem{Baxter} J.~R.~Baxter and  R.~V.~Chacon, \textit{The equivalence of diffussions on
networks to Brownian motion}, in  R.~Beals, A.~Beck, A.~Bellow, and
A.~Hajian (eds.), \textit{Conference on Modern Analysis and Probability},
Contemp. Math., Vol.~26,  Amer. Math. Soc., Providence, RI, 1984. p.~33 --
48.

\bibitem{Cardanobile} S.~Cardanobile and D.~Mugnolo, \textit{Analysis of a FitzHugh-Nagumo-Rall model of a neuronal network}, Math. Methods Appl. Sci. \textbf{30} (2007), 2281 -- 2308.

\bibitem{Carlson} R.~Carlson, \textit{Linear network models related to blood flow},
in G.~Berkolaiko, R.~Carlson,
S.~A.~Fulling, and P.~Kuchment (eds.), \textit{Quantum Graphs and Their
Applications}, Contemp. Math. Vol.~415, Amer. Math. Soc., Providence, RI,
2006. pp.~65 -- 80.

\bibitem{Chavel} I.~Chavel, \textit{Eigenvalues in Riemannian geometry},
    Pure and Applied Mathematics, Vol.~115. Academic Press, Orlando, 1984.

\bibitem{Chung} F.~R.~K.~Chung,
\textit{Spectral Graph Theory}, Regional Conference Series in Mathematics. Vol.~92. Providence, Amer. Math. Soc., Providence, RI, 1997.

\bibitem{Verdiere} Y.~Colin de Verdi\`{e}re,
\textit{Spectres de graphes}, Cours Sp\'{e}cialis\'{e}s, Vol.~4. Soci\'{e}t\'{e}
Math\'{e}matique de France, Paris, 1998.

\bibitem{Engel:Nagel} K.-J.~Engel and R.~Nagel, \textit{One-Parameter Semigroups for Linear
Evolution Equations}, Springer, New York, 2000.

\bibitem{Exner:Seba} P.~Exner and P.~{\v{S}}eba,
\textit{Free quantum motion on a branching graph}, Rep. Math. Phys.
\textbf{28} (1989), 7 -- 26.

\bibitem{Exner:Turek} P.~Exner and O.~Turek, \textit{Approximations of permutation-symmetric
vertex couplings in quantum graphs}, in G.~Berkolaiko, R.~Carlson,
S.~A.~Fulling, and P.~Kuchment (eds.), \textit{Quantum Graphs and Their
Applications}, Contemp. Math. Vol.~415, Amer. Math. Soc., Providence, RI,
2006. pp.~109 -- 120.

\bibitem{Mugnolo:4} M.~K.~Fijavz, D.~Mugnolo, and E.~Sikolya,
\textit{Variational and semigroup methods for waves and diffusion
in networks}, Appl. Math. Optimization \textbf{55} (2007), 219 --
240.

\bibitem{Freidlin:1}  M.~I.~Freidlin and A.~D.~Wentzell,
\textit{Diffusion processes on graphs and the averaging principle},
Ann. Probab.  \textbf{21} (1993),  2215 -- 2245.

\bibitem{Freidlin:2} M.~Freidlin and S.-J.~Sheu,
\textit{Diffusion processes on graphs: stochastic differential equations,
large deviation principle}, Probab. Theory Related Fields \textbf{116}
(2000), 181 -- 220.

\bibitem{Fulling} S.~A.~Fulling, P.~Kuchment, and J.~H.~Wilson, \textit{Index theorems for quantum graphs}, J. Phys. A: Math. Theor. \textbf{40} (2007), 14165 -- 14180.

\bibitem{Gaveau:1} B.~Gaveau, M.~Okada, and T.~Okada,
\textit{Explicit heat kernels on graphs and spectral analysis}, in
J.~E.~Fornaess (ed.), \textit{Several Complex Variables}, (Proceedings of
the Mittag-Leffler Institute, Stockholm, 1987-88), Princeton Math. Notes
Vol.~38, Princeton University Press, 1993. pp.~364 -- 388.

\bibitem{Gaveau:2} B.~Gaveau and M.~Okada,
\textit{Differential forms and heat diffusion on one-dimensional singular
varieties}, Bull. Sci. Math., II. S\'{e}r., \textbf{115} (1991), 61 -- 80.

\bibitem{Gilkey} P.~B.~Gilkey, \textit{Invariance Theory, the Heat Equation and the
    Atiyah-Singer Index Theorem}, Mathematics Lecture Series Vol.~11, Publish
    or Perish, Wilmington (1984).

\bibitem{Gilkey:2} P.~B.~Gilkey, \textit{Asymptotic formulae in spectral geometry},
Studies in Advanced Mathematics, Boca Raton, Chapman \& Hall/CRC, 2004.

\bibitem{Kochubej} A.~N.~Ko\v{c}ube{\u\i}, \textit{Extensions of symmetric operators and
    of symmetric binary relations}, Math. Notes \textbf{17} (1975), 25 -- 28.

\bibitem{Kato} T.~Kato, \textit{Perturbation Theory for Linear
Operators}, Springer-Verlag, Berlin, 1966.

\bibitem{KS1} V.~Kostrykin and R.~Schrader, \textit{Kirchhoff's rule
    for quantum wires}, J. Phys. A: Math. Gen. \textbf{32} (1999), 595 --
    630.

\bibitem{KS3} V.~Kostrykin and R.~Schrader, \textit{Kirchhoff's rule
    for quantum wires II: The inverse problem with possible
    applications to quantum computers}, Fortschr. Phys. \textbf{48} (2000),
    703 -- 716.

\bibitem{KS4} V.~Kostrykin and R.~Schrader, \textit{Quantum wires with magnetic fluxes},
    Comm. Math. Phys. \textbf{237} (2003), 161 -- 179.

\bibitem{KS8} V.~Kostrykin and R.~Schrader,
\textit{The inverse scattering problem for metric graphs and the traveling
salesman problem}, preprint \texttt{arXiv:math-ph/0603010} (2006).

\bibitem{KS9} V.~Kostrykin and R.~Schrader,
\textit{Laplacians on metric graphs: Eigenvalues, resolvents and
semigroups}, in G.~Berkolaiko, R.~Carlson, S.~A.~Fulling, and P.~Kuchment
(eds.), \textit{Quantum Graphs and Their Applications}, Contemp. Math.
Vol.~415, Amer. Math. Soc., Providence, RI, 2006. pp.~201 -- 225.

\bibitem{KPS1} V.~Kostrykin, J.~Potthoff, and R.~Schrader,
\textit{Heat kernels on metric graphs and a trace formula}, in F.~Germinet and P.~D.~Hislop (eds.), \textit{Adventures in Mathematical Physics}, Contemp. Math. Vol.~447, Amer. Math. Soc., Providence, RI, 2007. pp.~175 -- 198.

\bibitem{KPS2} V.~Kostrykin, J.~Potthoff, and R.~Schrader,
\textit{Brownian motion on metric graphs}, in preparation.

\bibitem{Kuchment:00} P.~Kuchment, \textit{Quantum graphs: I. Some basic structures},
Waves Random Media \textbf{14} (2004), S107 -- S128.

\bibitem{Kuchment:Post} P.~Kuchment and O.~Post, \textit{On the spectra of carbon nano-structures}, Commun. Math. Phys. \textbf{275} (2007), 805 -- 826.

\bibitem{Kurasov} P.~Kurasov, \textit{Graph Laplacians and topology}, Ark. Mat. (to appear).

\bibitem{Lumer} G.~Lumer, \textit{Connecting of local operators and evolution equations on networks}, in C.~Berg, G.~Forst, and B.~Fuglede (eds.), \textit{Potential theory} (Proceedings of a colloquium held in Copenhagen, May 14-18, 1979), Lect. Notes Math. Vol.~787. Berlin, Springer, 1980. pp.~219 -- 234.

\bibitem{Lumer:2} G.~Lumer, \textit{\'{E}quations de diffusion g\'{e}n\'{e}rales sur des r\'{e}seaux infinis}, in F.~Hirsch and G.~Mokobodzki (eds.), \textit{Seminaire de th\'{e}orie du potentiel},
    Lecture Notes in Math., Vol.~1061, Springer, Berlin, 1984. pp.~230 -- 243.

\bibitem{Mugnolo:2} D.~Mugnolo and S.~Romanelli,
\textit{Dynamic and generalized Wentzell node conditions for
network equations}, Math. Methods Appl. Sci. \textbf{30} (2007),
681 -- 706.

\bibitem{Mugnolo:3} D.~Mugnolo, \textit{Gaussian estimates for a heat equation on a
network}, Networks and Heterogeneous Media \textbf{2} (2007), 55 -- 79.

\bibitem{Nicaise:0} S.~Nicaise, \textit{Some results on spectral theory over networks, applied to nerve impulse transmission}, in C.~Brezinski, A.~Draux, A.~P.~Magnus, P.~Marino, and A.~Ronveaux (eds.), \textit{Polyn\^{o}mes orthogonaux et applications} (Proceedings of the Laguerre Symposium held at Bar-le-Duc, France, October 15-18, 1984).
Lect. Notes Math. Vol.~1171, Berlin, Springer, 1985. pp.~532 -- 541.

\bibitem{Okada} T.~Okada, \textit{Asymptotic behavior of skew conditional heat kernels
on graph networks}, Can. J. Math. \textbf{45} (1993), 863 -- 878.

\bibitem{Pankrashkin} K.~Pankrashkin, \textit{Spectra of Schr\"{o}dinger operators on equilateral quantum graphs}, Lett. Math. Phys. \textbf{77} (2006), 139 -- 154.

\bibitem{Phillips} R.~S.~Phillips, \textit{Dissipative operators and
hyperbolic systems of partial differential equations}, Trans. Amer. Math.
Soc. \textbf{90} (1959), 193 -- 254.

\bibitem{RS} M.~Reed and B.~Simon, \textit{Methods of Modern Mathematical Physics.
IV: Analysis of Operators}, Academic Press, New York, 1978.

\bibitem{Roth:1} J.-P.~Roth, \textit{Spectre du laplacien sur un graphe},
    C. R. Acad. Sci. Paris, S\'{e}r. I Math. \textbf{296} (1983), 793 -- 795.

\bibitem{Roth:2} J.-P.~Roth, \textit{Le spectre du laplacien sur un graphe},
    in G.~Mokobodzki and D.~Pinchon (eds.), \textit{Th\'{e}orie du potentiel},
    (Proceedings of the Colloque Jacques Deny, Orsay, June 20-23, 1983), Lecture
    Notes in Math., Vol.~1096, Springer, Berlin, 1984. pp.~521 -- 539.

\bibitem{Nagy:Foias} B.~Sz.-Nagy and C.~Foias, \textit{Harmonic Analysis of
Operators on Hilbert Spaces}, Amsterdam, North-Holland Publishing Company,
1970.

\bibitem{Tsekanovskii:2} E.~R.~Tsekanovskii, \textit{Non-self-adjoint
accretive extensions of positive operators and theorems of
Friedrichs-Krein-Phillips}, Funct. Anal. Appl. \textbf{14} (1980), 156 --
157.

\bibitem{Tsekanovskii} E.~R.~Tsekanovskii, \textit{Accretive extensions and
problems on the Stieltjes operator-valued functions relations}, in T.~Ando
and I.~Gohberg (eds.), \textit{Operator Theory and Complex Analysis}, Oper.
Theory Adv. Appl. Vol.~59, Basel, Birkh\"{a}user, 1992. pp.~328 -- 347.

\bibitem{Walsh} J.~B.~Walsh, \textit{A diffusion with a discontinuous
    local time}, Ast\'{e}risque \textbf{52-53} (1978), 37 -- 45.

\end{thebibliography}
\end{document}